%% file: indepnoout.tex
\def\hldest#1#2#3{}
\input fontmac
\input mathmac

\input epsf
\input eplain

\def\peel{\rho}
\def\leaf{\lambda}
\def\lp{\ell^+}
\def\lm{\ell^-}
\def\lss{\ell^{**}}
\def\rplus{r^+}
\def\rminus{r^-}
\def\rpe{r^{+\rm even}}
\def\rpo{r^{+\rm odd}}
\def\rme{r^{-\rm even}}
\def\rmo{r^{-\rm odd}}

\def\TV{{\rm TV}}
\def\l{\ell}
\def\Geo{{\rm Geo}}
\def\Uni{{\rm Uni}}
\def\Pos{{\rm Poi}}
\def\Bin{{\rm Bin}}
\def\TV{{\rm TV}}

\def\bref#1{[#1]}

\catcode`@=11
\dimendef\pl@left=0 \dimendef\pl@down=1
\dimendef\pl@right=2 \dimendef\pl@temp=3
\def\sob#1#2#3#4#5{
\setbox0\hbox{#1}\setbox1\hbox{$_\lhook$}\setbox2\hbox{p}%
\pl@right=#2\wd0 \advance\pl@right by-#3\wd1
\pl@down=#5\ht1 \advance\pl@down by-#4\ht0
\pl@left=\pl@right \advance\pl@left by\wd1
\pl@temp=-\pl@down \advance\pl@temp by\dp2 \dp1=\pl@temp
\leavevmode\kern\pl@right\lower\pl@down\box1\kern-\pl@left #1}
\def\eob{\sob e{.50}{.35}{0}{.93}}
\catcode`@=12

\font\mathbb=msbm10

\def\pr{\hbox{\mathbb P}}
\def\ex{\hbox{\mathbb E}}
\def\var{\hbox{\mathbb V}}

\enablehyperlinks

\def\ref#1{\special{ps:[/pdfm { /big_fat_array exch def big_fat_array 1 get 0
0 put big_fat_array 1 get 1 0 put big_fat_array 1 get 2 0 put big_fat_array pdfmnew } def}%
[\hlstart{name}{}{bib#1}#1\hlend]}

\baselineskip=13pt

\firstnopagenum

\maketitle{On the peel number and the leaf-height of Galton--Watson trees}{}
{Luc Devroye, Marcel K. Goh, {\rm and} Rosie Y. Zhao}{\sl School of Computer Science, McGill University}

\floattext5 \ninebf Abstract. \ninepoint\baselineskip=11.5pt
We study several parameters of a random Bienaym\'e--Galton--Watson tree $T_n$ of size $n$ defined in
terms of an offspring distribution $\xi$ with mean $1$ and nonzero finite variance $\sigma^2$. Let
$f(s)=\ex\{s^\xi\}$ be the generating function of the random variable $\xi$. We show that
the independence number is in probability asymptotic to $qn$, where $q$ is the unique solution to $q = f(1-q)$.
One of the many algorithms for finding the largest independent set of nodes uses a notion of repeated peeling
away of all leaves and their parents. The number of rounds of peeling is shown to be in probability asymptotic
to $\log n \big/ \log\big(1/f'(1-q)\big)$. Finally, we study a related parameter which we call the leaf-height.
Also sometimes called the protection number, this is
the maximal shortest path length between any node and a leaf in its subtree. If $p_1 = \pr\{\xi=1\}>0$,
then we show that the maximum leaf-height over all nodes in
$T_n$ is in probability asymptotic to $\log n/\log(1/p_1)$. If $p_1 = 0$
and $\kappa$ is the first integer $i>1$ with $\pr\{\xi=i\}>0$, then the leaf-height is in probability asymptotic
to $\log_\kappa\log n$.
\smallskip
\noindent{\bf Keywords.} Independence number, Bienaym\'e--Galton--Watson trees, protection number.

\advsect Introduction

{\tensc The independence number} is a fundamental graph invariant that arises often in computational complexity
theory and the analysis of algorithms. In a graph $G = (V,E)$, a subset $S\subseteq V$ of vertices is said to be
an {\it independent set\/} if no two elements of $S$ are adjacent. The dual notion is that of a {\it vertex cover},
namely a subset $C\subseteq V$ such that every edge in $G$ has an endpoint in $C$. The {\it independence number}
$I(G)$ of $G$ is defined to be the size of the largest independent set in $G$. In this paper, we concern ourselves
with the case $G=T$, a random tree in the Bienaym\'e--Galton--Watson model.
In recent years, analysis of the independence number
of trees has been carried out for various other random models. C.~Banderier, M.~Kuba, and
A.~Panholzer studied various families of simply-generated trees~\bref{4},
and a recent paper of M.~Fuchs, C.~Holmgren, D.~Mitsche, and R.~Neininger considers random binary search
trees as well as random recursive trees~\bref{13}.

Because every tree $T$ is bipartite, the independence number $I(T)$ is always at least $|T|/2$ (we
take the larger element of the bipartition).
Recall that a vertex set $S$ is a {\it vertex cover} of $T$ if every edge of $T$ intersects a vertex
in $S$. Letting $V(T)$ denote the size of a minimum-cardinality
vertex cover, we have the formula $n=V(T) + I(T)$.
In a tree, there always exists a maximum-cardinality
independent set that includes all of the leaves, and the following
algorithm, which will be the starting point of our discussion,
uses this fact to find an independent set of maximum size.
Note that this is only one of many possible algorithms that accomplishes this task.

\algbegin Algorithm I (Independent set). Given a directed tree $T$,
this algorithm computes a maximum-cardinality independent set $A$ of vertices.
\algstep I1. [Initialize.] Set $A \gets \emptyset$.
\algstep I2. [Compute leaves and parents.] Let $L(T)$ be the set of leaves of $T$, that is, the set of vertices
with out-degree $0$. Let $P(T)$ be the set of parents of nodes in $L(T)$.
\algstep I3. [Update.] Set $A \gets A\cup L(T)$ and $T\gets T\setminus L(T) \setminus P(T).$ (At this
stage, $T$ may now be a forest.)
\algstep I4. [Loop?] If $T = \emptyset$, halt and output $A$; otherwise, return to step I2.\slug

Algorithm I repeatedly peels away leaves and their parents to arrive at what we shall call the {\it layered
independent set}. We refer to $L(T)$ as layer~0, to $P(T)$ as layer~1, to $L\big(T \setminus L(T)
\setminus P(T)\big)$ as layer~2, and so on. In this manner, each node $u$ gets assigned a {\it peel number}
$\peel(u)$, the layer number of the set to which it belongs. The peel number $\peel(T)$
of a tree $T$ is the peel number of the root of $T$. We also let $m(T)$ denote the maximum of
the peel numbers of vertices in $T$; this quantity is twice the number of loops that Algorithm I undergoes before
termination, rounded up.
Note that all the peel numbers can be computed by postorder traversal of the tree
in time $O\big(|T|\big)$, and then the layered independent set is simply the collection of all nodes with even
peel number.

\newcount\figdefs
\figdefs=\figcount
A quantity related to the peel number is the {\it leaf-height} $\leaf(u)$
of a node $u\in T$. It is the length of the path to the nearest
leaf in the (fringe) subtree rooted at $u$. The leaf-height $\leaf(T)$
of a tree $T$ is the maximal leaf-height of any node in $T$. The fact that $\rho(u) = k$ implies that there is
a leaf at depth $k$ from the root, so $\leaf(u) \leq \rho(u)$ for all nodes $u$ in a tree. It is also easily
seen that for any tree $T$, $\leaf(T) \leq \rho(T)\leq m(T)$. A small example is given in Fig.~{\the\figdefs};
note that for nodes with few children or small subtrees, the two quantities are quite similar. One corollary of
our main results is that under certain conditions, this phenomenon persists as $n$ gets large, that is,
the peel number and leaf-height have the same order of asymptotic growth.
\midinsert
\vskip5pt
$$\epsfbox{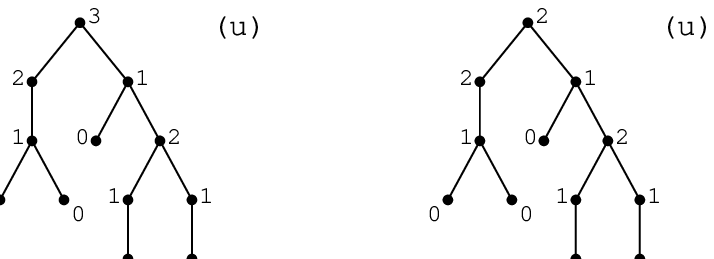}$$
\caption{Peel numbers and leaf-heights of nodes in a unary-binary tree.}
\endinsert
The leaf-height goes by the name {\it protection number} in the literature and has enjoyed some recent attention.
With this usage, a node whose minimal distance from any leaf is $k$ is called {\it $k$-protected} and a
$2$-protected node is often simply said to be {\it protected}.
In this paper we say that a node has leaf-height $k$, which we believe is more illustrative than saying it
is $k$-protected.
The number of nodes with leaf-height $\geq 2$ was examined by G.-S.~Cheon and L.~W.~Shapiro for
planted plane trees, Motzkin trees, full binary trees, Catalan trees, and ternary trees~\bref{8};
by T.~Mansour~\bref{23} for $k$-ary trees; by R.~R.~X. Du and and H.~Prodinger for digital search
trees~\bref{10}; by H.~M.~Mahmoud
and M.~D.~Ward for binary search trees~\bref{21} and
for random recursive trees~\bref{22}; and by L.~Devroye and
S.~Janson, who considered simply generated trees and also unified some earlier results regarding binary
search trees and random recursive trees. Nodes with leaf-height $>2$ were studied in binary search trees
by M.~B\'ona~\bref{5} and in planted plane trees by K.~Copenhaver~\bref{9}.
In the setting of simply generated trees and P\'olya trees,
the leaf-height of the root as well as the leaf height of a node chosen uniformly at random
was studied in~\bref{14}.

The main results of this paper characterize the asymptotic behaviour
of the independence number $I_n = I(T_n)$, the maximum peel number
$M_n = m(T_n)$, and the maximum leaf-height
$L_n = \lambda(T_n)$ for a Bienaym\'e--Galton-Watson tree $T_n$, which is conditioned on having $n$ nodes.
We also include
distributional properties of closely related statistics, such as the peel number and leaf-height
of the root of an unconditional Bienaym\'e--Galton-Watson tree, as well as
the leaf-height $L_n'$ of the root of and the leaf-height $L_n''$ of a node chosen uniformly at random in
a conditional Bienaym\'e--Galton--Watson tree.

\medskip\boldlabel The Bienaym\'e--Galton--Watson model.
For a nonnegative integer-valued random variable $\xi$, a {\it Bienaym\'e--Galton--Watson tree}
is a random tree in which every node has $i$ children independently with
probability $p_i = \pr\{\xi=i\}$. The random variable $\xi$ is called the {\it offspring distribution} of
the tree; we only consider distributions with mean $\ex\{\xi\}
=1$ and variance $\var\{\xi\} = \sigma^2 \in (0,\infty)$ (standard references include~\bref{3}
and~\bref{20}). Let $T_n$ denote the tree
$T$, conditioned on having $n$ nodes. Note that many important
simply generated families of trees can be characterized
by a conditional Bienaym\'e--Galton--Watson tree with a certain distribution~\bref{15}
Strictly speaking, $T$ is a graph $(V,E)$,
but we will abuse notation and write $v\in T$ to indicate that $v$ is in the vertex set of $T$.

\advsect The independence number

We begin by studying unconditional Bienaym\'e--Galton--Watson trees.
Recall that the generating function $f(s)$ of an offspring
distribution $\xi$ is the infinite series $\ex\{s^\xi\}$, which converges absolutely when $0\leq s\leq 1$. We can
thus differentiate to obtain $f'(s) = \ex\{\xi s^{\xi-1}\}$. A quantity that will play a key role in our story is
$q$, the unique solution in $(0,1)$ of $q = f(1-q)$.
\midinsert
\vskip5pt
$$\epsfbox{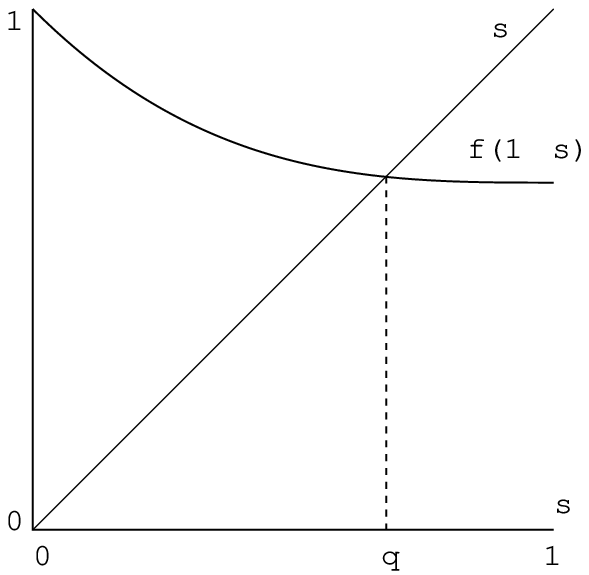}$$
\caption{The parameter $q$ satisfying $q = f(1-q)$.}
\endinsert
\newcount\indepprob
\indepprob=\thmcount
\proclaim Lemma \advthm. Let $\xi$ be an offspring distribution with $0<\ex\{\xi\} \leq 1$ and let
$f(s)=\ex\{s^\xi\}$. The probability that the root of a
Bienaym\'e--Galton--Watson tree $T$ with this distribution belongs
to the layered independent set is $q$, which belongs to the interval $(1/2,1)$.

\proof
Note that $q$ is the probability that all the children of the root are not in the layered independent set.
By the recursive definition of a Bienaym\'e--Galton--Watson tree, we have
$$q = \sum_{i\geq 0} p_i(1-q)^i = f(1-q),\adveq$$
and the Banach fixed-point theorem guarantees the uniqueness of the solution to $s = f(1-s)$ in the compact interval
$[0,1]$. Of course, $q$ cannot be $1$ since $\pr\{\xi = 0\}\neq 0$. The fact that $f(s) > s$ for all
$s\in (0,1)$ implies that $q = f(1-q) > 1-q$, hence $q>1/2$.\slug
\medskip
Lemma~{\the\indepprob} is essentially known (see, e.g., Banderier, Kuba, and Panholzer~\bref{4}).
\medskip
\boldlabel Examples. There is a well-known connection between certain families of trees and conditioned
Bienaym\'e--Galton--Watson trees. In each of the following cases, sampling a conditional
Bienaym\'e--Galton--Watson tree $T_n$
with the given distribution
is equivalent to uniformly sampling a tree of size $n$ from the respective tree family.
\medskip
\item{i)} In {\it Flajolet's $t$-ary tree}, every node is either a leaf or has $t$ children. This corresponds
to the distribution with $p_0 = 1-1/t$ and $p_t = 1/t$, so we can compute $q$ numerically by finding
the unique solution to the equation
$$q = 1-{1\over t} + {(1-q)^t\over t}.\adveq$$
in the interval $(1/2,1)$. In the case $t=2$ of full binary trees, we find that $q = 2-\sqrt 2 \approx 0.585786$,
and since the $(1-q)^t/t$ term is very small for larger values of $t$, $q$ is approximately $1-1/t$ for large $t$.
\smallskip
\item{ii)} To obtain a random rooted {\it Cayley tree},
we set $p_i = (i!e)^{-1}$ for all $i\geq 0$. Since $f(s)=e^{s-1}$,
we have $qe^q = 1$, which we can invert in terms of the {\it Lambert $W$ function}. Concretely, we have
$$q = W(1) = \bigg(\int_{-\infty}^\infty {dt\over (e^t-t)^2 + \pi^2}\bigg)^{-1} - 1 \approx 0.567143,\adveq$$
which is also known as the {\it omega constant}.
\smallskip
\item{iii)} {\it Planted plane trees} correspond to the distribution $p_i = 1/2^{i+1}$ for $i\geq 0$. In this case,
$f(s) = 1/(2-s)$, yielding the equation $q^2+q-1 = 0$, whose solution in the correct range is $q=1/\varphi \approx
0.618034$. (The golden ratio $\varphi = 1.618034$ is the more famous solution to this quadratic equation).
\smallskip
\item{iv)} {\it Motzkin trees}, also known as unary-binary trees, are trees in which every non-leaf node has
either one tree or two children. This corresponds to the distribution $p_0 = p_1 = p_2 = 1/3$ and $p_i = 0$ for
all $i\geq 3$. So we have $q = \big(1 + (1-q) + (1-q)^2\big)/3$ and we have $q = 3-\sqrt 6 \approx 0.550510$.
\smallskip
\item{v)} A {\it binomial tree} of order $d$ can be thought of as a tree in which every node has $d$ ``slots''
for its children, some of which may be filled. Thus a node can have $r$ children in ${d\choose r}$ different ways,
for $0\leq r\leq d$. This corresponds, fittingly, to a binomial offspring distribution, where
$$p_i = {d\choose i}\bigg({1\over d}\bigg)^i\bigg(1-{1\over d}\bigg)^{d-i},\adveq$$
for $0\leq i\leq d$, and $p_i = 0$ otherwise. For this distribution, we have $f(s) = (s/n + 1- 1/n)^n$, meaning
that
$$q = \bigg(1-{1\over d} + {1-q\over d}\bigg)^d = \bigg(1-{q\over d}\bigg)^d.\adveq$$
For large $d$, this tends to the omega constant. An important case is $d=2$, which produces a random {\it Catalan
tree}; it can be readily computed that $q = 4-2\sqrt 3\approx 0.535898$ for these trees.
\medskip\goodbreak

The following theorem shows the link between $q$ and the size of the largest independent set in a conditioned
Bienaym\'e--Galton--Watson tree.

\proclaim Theorem \advthm. Let $\xi$ be an offspring distribution with $\ex\{\xi\}= 1$ and let
$f(s)=\ex\{s^\xi\}$. The independence number $I_n = I(T_n)$ of a Bienaym\'e--Galton--Watson tree, conditioned on having
$n$ nodes, satisfies
$${I_n\over n} \to q$$
in probability as $n\to \infty$, where $q$ is the unique solution in $(1/2,1)$ of the equation $q=f(1-q)$.

\proof For a vertex $u$, we let $\Gamma_u$ denote the set of children of $u$ and let
$$g(u) = \cases{1,& if the peel number of $u$ is even;\cr 0,&otherwise}.\adveq$$
Note that the recursive function
$$G(u) = g(u) + \sum_{v\in \Gamma_u} G(v)\adveq$$
is exactly the independence number of the subtree rooted at $u$. Since $g$ is bounded, we can apply a result of
S.~Janson (\bref{16}, Theorem 1.3) to find that for a conditional
Bienaym\'e--Galton--Watson tree $T_n$ with root $u$,
$${I_n\over n} = {G(u)\over n}  \to \ex\big\{g(u)\big\} = q\adveq$$
in probability as $n\to\infty$.\slug

Note that examples (ii), (iii), and the Catalan case agree with explicit computations given
in~\bref{4}. For simply generated trees, that paper, which uses singularity analysis,
derives the constant $q$ in a different manner, proves the
stronger statement $\ex\bigl\{I(T_n)\bigr\} = qn + O(1)$, and also gives a formula for the variance in terms
of the degree-weight generating function. In particular, they show there exists a constant $\nu$ depending
on the family of trees such that the variance is $\nu n + O(1)$.

\newcount\spathsect
\spathsect=\sectcount
\advsect Minimum-size {\mathbold s}-path vertex covers

This section represents a brief digression, and will not be related to our remaining results, though it discusses
the natural generalization of Algorithm I and is related to the open problem we give at the end of the paper.
As mentioned in the introduction, the size $V(T)$ of the minimal vertex cover of a $T$ with $n$ nodes
has size $n - I(T)$, where $I(T)$ is the independence number. In particular, Algorithm I outputs a
minimum-cardinality
vertex cover alongside the maximum-cardinality
independent set; it is the set of all nodes with odd peel number. We now
tackle a more general notion of vertex covers. For an integer $s\geq 2$, an {\it $s$-path vertex cover} of a rooted
tree $T$ is a subset $C$ of vertices such that any path of length $s-1$ in the tree contains a vertex in $C$.
Thus the common-or-garden vertex cover corresponds to $s=2$.
(The off-by-one quirk in the definition goes away
if we measure a path not by its length, but instead by its {\it order}, that is, the number of vertices it
contains.) Note that only directed paths are considered; so two children of the same node are {\it not}
connected by a path of length 2.

One might be tempted to generalize our earlier observations by claiming that the set of nodes with peel number
congruent to $s-1$ modulo $s$ is a minimal $s$-path vertex cover. This is not true! Consider a tree
in which the root has two children, and one of the children has itself one child. Then no node has peel number
equal to $2$, but of course, the minimal $3$-path vertex cover consists of the root. Towards a correct
generalization, consider the fact that if in every loop of Algorithm I, we removed all subtrees of height
exactly 1, then the roots of these removed subtrees are precisely the vertices with odd peel number. Thus
we arrive at an algorithm for computing a minimal $s$-path vertex cover.

\algbegin Algorithm P (Compute $s$-path vertex cover). Given $s\geq 2$, and a rooted tree $T$, this algorithm
computes a minimal $s$-path vertex cover $C$.
\algstep P1. [Initialize.] Set $C\gets \emptyset$.
\algstep P2. [Done?] If there are no subtrees with height exactly $s-1$, we output $C$ and terminate.
\algstep P3. [Prune a subtree.] Let $v$ be a node in $T$ such that the subtree $T_v$ rooted at $v$ has
height exactly $s-1$. We set $C\gets C\cup\{v\}$ and set $T\gets T\setminus T_v$. Return to step P2.\slug

Note that if the original tree had height less than $s-1$, the algorithm outputs the empty set, which
is a valid cover, since there are no paths of length $s-1$ in the tree. The fact that this algorithm actually
does output a minimum-size vertex cover
is proved in~\bref{6}, and it is also remarked that the algorithm can be made to run in
$O\big(|T|\big)$ time.

Let $V_s(T)$ denote the size of the minimum $s$-path vertex cover of a
Bienaym\'e--Galton--Watson tree $T$. To determine this
random quantity, we will have to determine the probability that a node is added to the set $C$ in Algorithm P.
We will say that a vertex $v\in T$ is ``marked'' if Algorithm P adds it to the cover $C$. The following
lemma gives necessary and sufficient conditions for the root of a tree to be marked.

\proclaim Lemma \advthm. The root $u$ of a tree is marked if and only if there exists a path of length $s-1$
from the root that contains no marked vertices (other than the root).

\proof Suppose that the root $u$ is marked. This means that in the final iteration of Algorithm P, after
all other marked nodes have been removed, the tree has height $s-1$. This means that some unmarked node $v$
is at depth $s-1$, and no node is marked on the path to this node. (This happens when $v$
is a leaf or all children of $v$ are marked, since if a child of $v$ is unmarked, then we have an unmarked path
of length $s$ in the tree and the algorithm would have to mark some node on this path before marking
the root.) Conversely, if such a path exists,
then Algorithm P will be in this state in the final iteration of the loop and will therefore mark the
root.\slug

This observation can be used to derive a functional equation for the probability that a node in an
unconditional Bienaym\'e--Galton--Watson tree is marked, as the following lemma shows.

\newcount\spath
\spath=\thmcount
\proclaim Lemma \advthm. Let $T$ be a Bienaym\'e--Galton--Watson tree with offspring distribution $\xi$ satisfying
$\ex\{\xi\}\leq 1$. Let $f(z) = \ex\{\xi^z\}$ be the generating function of the distribution and let
$$g(z, q) = 1-f\big(q+(1-q)z\big)\adveq$$
The probability $q_s$ that the root of the tree is in the minimum $s$-path vertex
cover produced by Algorithm P satisfies
$$q_s = g(g(\cdots(g(0,q_s)\cdots,q_s),q_s),q_s),\adveq$$
where the function $g$ is iterated $s-1$ times.

\proof For $1\leq j<s$, let $E_j$ be the event that in an unconditional Bienaym\'e--Galton--Watson tree,
there is a path of length $j$ from the root that contains no marked vertices (except possibly the root).
Thus $q_s = \pr\{E_{s-1}\}$. Restating things slightly, $E_j$ is the probability that there exists an unmarked
child $v$ of the root in whose subtree $E_{j-1}$ is true. If the degree of the root is $i$,
then the probability that all the children of the root fail to have this property is
$\big(q_s + (1-q_s)\pr\{E_{j-1}\}\big)^i$, so for $1<j<s$,
$$\pr\{E_j\} = \sum_{i\geq 0} p_i\big(1-(q_s+(1-q_s)\pr\{E_{j-1}\})^i\big)
= 1-f\big(q_s+(1-q_s)\pr\{E_{j-1}\}\big).\adveq$$
Note that $\pr\{E_1\}$ is simply the probability $1-f(q_s) = g(0, q_s)$ that one of the children
of the root is unmarked, so unravelling the above equation proves the lemma.\slug

Note that when $s=2$, $q_s = 1-q$, where $q$ is the solution to $z=f(1-z)$ we studied earlier.
By a recursive computation analogous to the one we performed for the independence number, we find that
if $V_s(T_n)$ denotes the minimum size of an $s$-path vertex cover of the
conditional Bienaym\'e--Galton--Watson tree
$T_n$, then as $n\to \infty$,
$${V_s(T_n)\over n}\to q_s,\adveq$$
in probability. The function $g$ given by Lemma~{\the\spath} is rather unwieldy, so we cannot hope to
find neat closed forms for the limit of $V_s(T_n)$ like we did for $I_n$ in many special cases. However,
we can, in principle, use $g$ to numerically approximate the $s$-path vertex cover number for arbitrary
distributions satisfying $\ex\{\xi\} \leq 1$.

\newcount\peelnumsect
\peelnumsect=\sectcount
\advsect Distribution of the peel number

Let $r_i$ denote the probability that the root of an unconditional
Bienaym\'e--Galton--Watson tree has peel number $i$. In this section,
we shall compute the distribution $(r_i)_{i\geq 0}$.
It will also be convenient to set $r_i = 0$ when $i$ is negative. We will establish the notation
$$\rplus_i= \sum_{j\geq i}r_j\qquad\hbox{and}\qquad \rminus_i = \sum_{j=0}^i r_j\,.\adveq$$
There will be some asymmetry for odd and even $i$, so let
us write $\rpo_i$ for the subsum of $\rplus_i$ consisting of odd terms and $\rpe_i$ for the subsum of $\rplus_i$
consisting of even terms. Defining $\rmo_i$ and $\rme_i$ similarly, we have, of course,
$\rpo_i+\rpe_i=\rplus_i$ and $\rmo_i + \rme_i = \rminus_i$.
\midinsert
\vskip5pt
$$\epsfbox{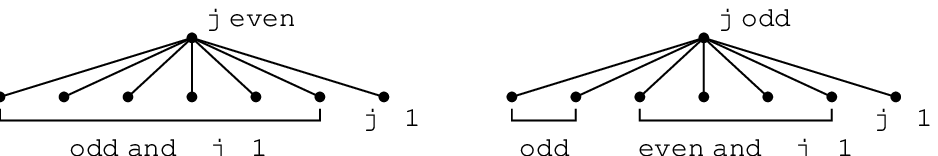}$$
\caption{Children of nodes with even and odd peel numbers.}
\endinsert
\noindent Clearly, $r_0 = p_0$. For even indices $j$, all children must have an odd peel number at most $j-1$
and at least one must
have peel number $j-1$. Thus, if $\xi$ is the number of children at the root, then for $i\geq 1$,
$$r_{2i} = \ex\big\{\big(\rmo_{2i-1}\big)^\xi-\big(\rmo_{2i-3}\big)^\xi\big\} = f\big(\rmo_{2i-1}) -
f\big(\rmo_{2i-3}\big).\adveq$$
For odd indices $j$, all the children of the root with even peel number must have peel number at least
$j-1$, and at least one must have peel number $j-1$. Since
$$\sum_{i\geq 0} r_{2i-1} = 1-q\qquad\hbox{and}\qquad \sum_{i\geq 0} r_{2i} = q,\adveq$$
we find that for $i\geq 1$,
$$r_{2i-1} = \ex\big\{\big(1-q+\rpe_{2i-2}\big)^\xi - \big(1-q+\rpe_{2i}\big)^\xi\big\}.\adveq$$
The following lemma describes $r_i$ for large $i$.

\proclaim Lemma \advthm. Let $r_i$ be the probability that the root of an unconditional Bienaym\'e--Galton--Watson
tree with offspring distribution $\xi\sim(p_i)_{i\geq 0}$ has peel number equal to $i$. As $i\to \infty$,
we have
$$r_i = f'(1-q)^{i+o(i)}.\adveq$$

\proof
In the even case, we have
$$\eqalign{
r_{2i} &= f\big(\rmo_{2i-1}\big) - f\big(\rmo_{2i-3}\big) \cr
&\sim r_{2i-1}\sum_{j\geq 0} jp_j \big(\rmo_{2i-3}\big)^{j-1}\cr
&= r_{2i-1} f'\big(\rmo_{2i-3}\big),\cr
}\adveq$$
which, since $\rmo_{2i-3}\to 1-q$, is asymptotic to $r_{2i-1}f'(1-q)$. Similarly, we have
$$r_{2i-1} = f\big(\rpe_{2i-2} + 1-q\big) - f\big(\rpe_{2i} + 1-q\big)
  \sim r_{2i-2}\sum_{j\geq 0} jp_j(1-q)^{j-1},\adveq$$
which is also asymptotic to $r_{2i-2}f'(1-q)$.\slug

If $N_i$ is the number of nodes in the $i$th layer for our algorithm, then Aldous's theorem~\bref{1}
implies that for every fixed $i$,
$${N_i\over n}\to r_i\adveq$$
in probability. The number of nodes in the layers decreases at the indicated rate, namely $f'(1-q)$. As
$q\in (1/2,1)$, we have
$$p_1 = f'(0) < f'(1-q) < f'\bigg({1\over 2}\bigg) \leq \ex\bigg\{{1\over 2^\xi}\bigg\}.\adveq$$

The next section will need the event that the maximum peel number in an unconditional tree occurs at the root. We have the following lemma.

\newcount\rootmax
\rootmax=\thmcount
\proclaim Lemma \advthm. Let $T$ be an unconditional Bienaym\'e--Galton--Watson tree with offspring distribution
$\xi$.
Let $R$ be the peel number of the root of such a tree and let $M$ be the maximum peel
number of any node in the tree.
Let $q$ be the solution to $q=f(1-q)$, where $f$ is the reproduction
generating function of this distribution. Then
$$\tau_i := \pr\{R = M = i\} = f'(1-q)^{i+o(i)}\adveq$$
as $i\to\infty$.

\proof The fact that $\tau_i \leq \pr\{R = i\} = r_i = f'(1-q)^{i+o(i)}$ means that we only have to worry
about finding a lower bound. To that end, consider the $\xi$ children of the root (each the root of
unconditional Bienaym\'e--Galton--Watson trees), with peel numbers
$R_1,\ldots,R_\xi$ and maximum peel numbers $M_1,\ldots, M_\xi$. We consider the odd and even cases separately.

When $i$ is odd, the event that $R=M=i$ is
implied by the event that there exists some $1\leq j\leq \xi$ with $R_j = M_j = i-1$ and for all $k\ne j$,
we have $R_j$ odd and $M_j\le i$. Therefore, when $i$ is odd, we have, by the inclusion-exclusion inequality,
$$\tau_i \ge
\ex\big\{ \xi\cdot \tau_{i-1} \pr\{R\ \hbox{odd},\, M\le i\}^{\xi-1}\big\}
-\ex\bigg\{ {\xi\choose 2} \cdot {\tau_{i-1}}^2 \pr\{R\ \hbox{odd},\, M\le i\}^{\xi-2}\bigg\}.\adveq$$
Note that $\pr\{R\ \hbox{odd},\,M> i\} = o(1)$ as $i\to\infty$, and so
$$\eqalign{
\tau_i &\ge  \ex\big\{ \xi\cdot \tau_{i-1} \big(1-q-o(1)\big)^{\xi-1}\big\}
  -\ex\bigg\{ {\xi\choose 2} \cdot {\tau_{i-1}}^2 \big(1-q-o(1)\big)^{\xi-2}\bigg\} \cr
&= \tau_{i-1} f'\big(1-q-o(1)\big) - {{\tau_{i-1}}^2\over 2} f''(1-q) \cr
&\ge \tau_{i-1} f'\big(1-q-o(1)\big) - {{\tau_{i-1}}^2 \sigma^2\over 2}.\cr
}\adveq$$
When $i$ is even, the event that $R=M=i$ is implied by the event that there exists some $1\leq j\leq \xi$
with $R_j=M_j =i-1$ and for all $k\ne j$, we have $R_j$ odd, $R_j\leq i-2$, and $M_j\leq i$. With another
application of the inclusion-exclusion inequality and by a similar argument as in the odd case, we have
$$\eqalign{
\tau_i &\ge \ex\big\{ \xi\cdot \tau_{i-1} \pr\{R\le i-2,\,M\le i,\, R\ \hbox{odd}\big\}^{\xi-1}\big\}
  - \ex\bigg\{ {\xi\choose 2} {\tau_{i-1}}^2 \pr\{R\le i-2,\,M\le i,\, R\ \hbox{odd}\}^{\xi-2}\bigg\} \cr
&\ge \tau_{i-1} f'\big(1-q-o(1)\big) - {{\tau_{i-1}}^2\sigma^2\over 2}.\cr
}\adveq$$
In both the odd and even cases, we see that $\tau_i\ge f'(1-q)^{i+o(i)}$, completing the proof.\slug

Using this result, we can give the following property of the distribution of the maximum peel number.

\proclaim Lemma \advthm. The maximum peel number $M$ in an unconditional Bienaym\'e--Galton--Watson tree satisfies
$$\pr\{M\geq i\} = f'(1-q)^{i/2 + o(i)},\adveq$$
as $i\to\infty$, where $f$ is the reproduction generating function.

\proof As before, let $R_1,\ldots,R_\xi$ denote the peel numbers of children of the root and let
$M_1,\ldots, M_\xi$ denote the maximum peel numbers in their respective subtrees.
Let $\mu_i = \pr\{M=i\}$, $\mu_i^- = \pr\{M\le i\}$, and $\mu_i^+=\pr\{M\ge i\}$. The event that
$M\ge i$ is implied by the event
\newcount\event
\event=\eqcount
$$\Big(\max_{1\le j\le \xi} M_j\ge i\Big)\ \hbox{or}
\ \Big(\max_{1\le j\le\xi}M_j<i\ \hbox{and}
\ \hbox{there is some} \ 1\le j\le\xi\ \hbox{with}\ R_j=M_j=i-1\Big).\adveq$$
Thus, letting $E_j$ be the event that $R_j=M_j=i-1$, we have
\newcount\oneandtwo
\oneandtwo=\eqcount
$$\pr\{M\ge i\} \ge \pr\Big\{\max_{1\le j\le\xi} M_j\ge i\Big\}
   + \pr\bigg\{ \max_{1\le j\le\xi} M_j<i,\, \bigcup_{j=1}^\xi E_j\bigg\}.\adveq$$
Note first that
$$\pr\Big\{\max_{1\le j\le\xi} M_j\ge i\Big\} = 1 - \ex\big\{(1-\mu_i^+)^\xi\big\} = 1-f(1-\mu_i^+).\adveq$$
By taking the Taylor series expansion of $f$ around $1$, we have
$$f(1-s) = f(1) - sf'(1)+{s^2\over 2} f''(\theta)\adveq$$
for some $1-s\le\theta\le 1$, so that
$$f(1-s) =1-s+{s^2\over 2}f''(\theta)\leq 1-s+{s^2\over 2}\sigma^2\adveq$$
and
$$\pr\Big\{\max_{1\le j\le\xi} M_j\ge i\Big\} \ge \mu_i^+ - {{{\mu_i^+}^2}\sigma^2\over 2}.\adveq$$
Next, by the union bound, we have
$$\eqalign{
\pr\bigg\{ \max_{1\le j\le\xi} M_j<i,\, \bigcup_{j=1}^\xi E_j\bigg\} &=
  \pr\bigg\{\bigcup_{j=1}^\xi E_j\bigg\} - \pr\bigg\{\bigcup_{j=1}^\xi E_j,\,
  \max_{1\leq j\leq \xi} M_j\ge i\bigg\} \cr
&= 1 - \ex\big\{ (1-\tau_{i-1})^\xi\big\} - \ex\big\{ \xi(\xi-1) \tau_{i-1}\mu_i^+\big\} \cr
&= 1 - f(1-\tau_{i-1}) - \sigma^2 \tau_{i-1}\mu_i^+ \cr
&\ge \tau_{i-1} - {{\tau_{i-1}}^2\sigma^2\over 2} - \sigma^2 \tau_{i-1}\mu_i^+.\cr
}\adveq$$
Collecting these bounds back into~\refeq{\the\oneandtwo}, we have
$$\mu_i^+ \ge \mu_i^+ + \tau_{i-1} - {{\mu_i^+}^2\sigma^2\over 2} - {{\tau_{i-1}}^2\sigma^2\over 2}
-\sigma^2 \tau_{i-1}\mu_i^+\adveq$$
and therefore
\newcount\maxpeeleqseven
\maxpeeleqseven=\eqcount
$${\sigma^2\over 2}{\mu_i^+}^2 \geq \tau_{i-1}(1-\sigma^2\mu_i^+) - {{\tau_{i-1}^2\sigma^2}\over 2}.\adveq$$
Let $\phi(x)$ be a decreasing function that is $o(1)$ as $x\to\infty$. We combine~\refeq{\the\maxpeeleqseven}
with Lemma~{\the\rootmax} to conclude that
$${\sigma^2\over 2}{\mu_i^+}^2 \ge \tau_{i-1}\big(1-\phi(i)\big) = f'(1-q)^{i+o(i)}.\adveq$$
To bound $\mu_i^+$ from above, we observe that since the event that $M\ge i$ is a subset of the
event~\refeq{\the\event}, we have
$$\eqalign{
\mu_i^+ &\le 1-f(1-\mu_i^+) + \ex\big\{ 1-(1-\tau_{i-1})^\xi\big\} \cr
&\le \mu_i^+ - {{\mu_i^+}^2\over 2}\big(\sigma^2+o(1)\big) + \tau_{i-1} - {{\tau_{i-1}}^2\over 2}
\big(\sigma^2+o(1)\big),\cr
}\adveq$$
and therefore
$${\sigma^2+o(1)\over 2}{\mu_i^+}^2 \leq \tau_{i-1} = f'(1-q)^{i+o(i)},\adveq$$
which is what we need.\slug

\newcount\peelnumasymp
\peelnumasymp=\sectcount
\advsect Asymptotics of the peel number

We are now ready to prove an asymptotic result for the peel number of $T_n$; that is,
the maximum peel number over all nodes in $T_n$. This is the number of rounds of peeling required
by Algorithm~I to calculate a maximum-cardinality independent set.
Our proof uses Kesten's
tree $T_\infty$, whose construction we shall briefly recall here (see~\bref{18}).
Fix an offspring distribution $\xi$
with $\ex\{\xi\} = 1$. Starting from the root, we attach $\zeta$
children, where $\pr\{\zeta = i\} = ip_i$ for $i\geq 0$. Now select a child uniformly at random and mark it.
We repeat the process at the marked child, while all other children become the roots of independent unconditional
Bienaym\'e--Galton--Watson trees with the ordinary offspring distribution $\xi$.
For a tree $t$, we let $\tau(t,k)$ denote
$t$, truncated to include only the first $k$ levels. Letting $\TV$ denote total variation distance, it
is well known (see~\bref{17} and~\bref{25}) that if $k=o(\sqrt n)$, then
$$\lim_{n\to\infty} \TV\big(\tau(T_\infty,k), \tau(T_n,k)\big) = 0.\adveq$$

\newcount\peelthm
\peelthm=\thmcount
\proclaim Theorem \advthm. Let $M_n$ be the maximum peel number
in $T_n$, a conditional Bienaym\'e--Galton--Watson tree on $n$ nodes with offspring distribution $\xi$. Then
$${M_n\over\log n}\to {1\over \log\big(1/f'(1-q)\big)}\adveq$$
in probability, where $f$ is the generating function of $\xi$.

\proof For any tree $t$, let $h(t)$ denote its
height and $m(t)$ its maximum peel number. For the lower bound, we employ Kesten's limit tree $T_\infty$.
Let $S_k$ denote the set of nodes of $T_\infty$ that are children of nodes on the spine of $\tau(T_\infty,k)$ (i.e. nodes that are marked in the construction of $T_\infty$).
Let
$$\alpha_n = \bigg\lfloor {\sqrt n\over \log^2 n}\bigg\rfloor\qquad\hbox{and}\qquad
\beta_n = \bigg\lfloor {\sqrt n\over \log n}\bigg\rfloor.$$
By the same result of~\bref{17} and~\bref{25} that we used before, we can find a coupling
of $\tau(T_n, \beta_n)$ and $\tau(T_\infty, \beta_n)$ such that
$$\pr\big\{ \tau(T_n,\beta_n) \ne \tau(T_\infty,\beta_n)\big\} = o(1).\adveq$$
For every node $u$ in $T_\infty$, let $T_u$ be the subtree of $T_\infty$ rooted at $u$. Let $M_n = m(T_n)$.
Letting $E_{ux}$
denote the event that $m(T_u)\le x$, we have
$$\pr\{M_n\le x\} \leq \pr\big\{ \tau(T_n,\beta_n) \ne \tau(T_\infty,\beta_n)\big\}
  +\pr\Big\{ \bigcap_{u\in S_{\alpha_n}} E_{ux}\Big\}
  +\pr\Big\{\max_{u\in S_{\alpha_n}} h(T_u) \ge \beta_n-\alpha_n\Big\}.\adveq$$
We already pointed out that $\pr\big\{ \tau(T_n,\beta_n) \ne \tau(T_\infty,\beta_n)\big\} = o(1)$; we bound
the other two terms by
\newcount\fourterms
\fourterms=\eqcount
$$\eqalign{
\pr\Big\{ \bigcap_{u\in S_{\alpha_n}} E_{ux} \Big\}
  +\pr\Big\{\max_{u\in S_{\alpha_n}} h(T_u) \ge \beta_n-\alpha_n\Big\}
  &\le \pr\bigg\{|S_{\alpha_n}| \leq {\sigma^2 \alpha_n\over 2}\bigg\}
   + \pr\bigg\{|S_{\alpha_n}| \geq {\sigma^2 3\alpha_n\over 2}\bigg\}\cr
   &\quad+ \pr\big\{ m(T)\leq x\big\}^{\sigma^2\alpha_n/2} + {2\sigma^2\over 2}\alpha_n \pr\big\{h(T)\geq \beta_n-
   \alpha_n\big\},\cr
}\adveq$$
where $T$ is an unconditional Bienaym\'e--Galton--Watson tree. Now, $S_{\alpha_n}/ (\sigma^2\alpha_n)\to 1$
in probability by the law of large numbers, as the expected number of children of any node on the spine of
$T_\infty$ is $\sigma^2+1$. So, the first two terms of~\refeq{\the\fourterms} tend to zero. Next, we see that
$$\pr\big\{m(T)\le x\big\}^{\sigma^2\alpha_n/2} = \big(1-\pr\big\{m(T)> x\big\}\big)
^{\sigma^2\alpha_n/2} \leq \exp\!\bigg(\!\!-\! f'(1-q)^{x/2 + o(x)}{\sigma^2\alpha_n \over 2}\bigg),\adveq$$
which tends to zero if $x = (1-\eps)\log n/\log\big(1/f'(1-q)\big)$. For the final term, we have, by
Kolmogorov's theorem (see, e.g.~\bref{20} or~\bref{2})
$${3\sigma^2\alpha_n\over 2} \pr\big\{h(T) > \beta_n-\alpha_n\big\} \sim
{3\sigma^2\over 2} \cdot {2\alpha_n\over \sigma^2(\beta_n-\alpha_n)}\sim {3\alpha_n\over \beta_n}
\sim {3\over \log n},\adveq$$
which goes to zero. We have shown that
$$\pr\bigg\{ M_n < (1-\eps) {\log n\over \log \big(1/f'(1-q)\big)}\bigg\}\to 0\adveq$$
for all $\eps>0$.
\midinsert
\vskip5pt
$$\epsfbox{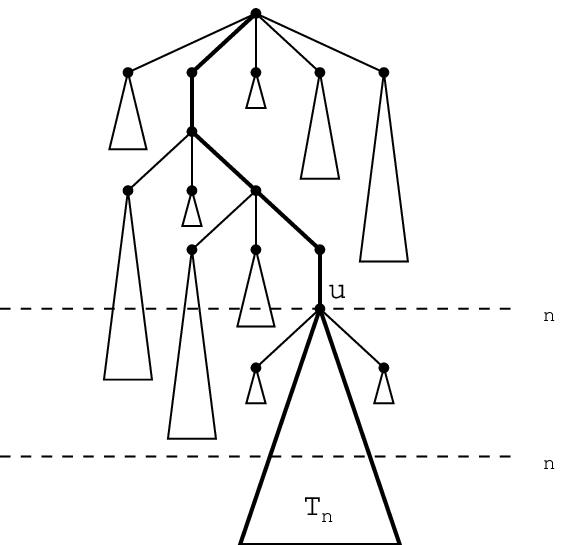}$$
\caption{The proof uses Kesten's infinite tree $T_\infty$ for both bounds.}
\endinsert
For the upper bound we will again work with $T_\infty$, truncated to level $\beta_n$, but also require some
further auxiliary definitions. Let $u^*$ denote the unique node on the spine of $T_\infty$ at distance $\alpha_n$ from
the root of $T_n$. Let $T_n^*$ be its subtree (in $T_n$, not $T_\infty$).
Let $S$ be the set of children of nodes on the spine
at distance $\leq \alpha_n$ from the root. We then define
$$M_n' = \max_{u\in S} \rho(u)\qquad\hbox{and}\qquad M_n'' = \max_{u\in S} m(T_u).$$
Next, we let $S^*$ denote the set of nodes $u$ on the spine with the property that all of $u$'s non-spine
children have an odd peel number. In particular, let $Y_n$ be the maximal number of {\it consecutive} nodes on
the spine that are in $S^*$. Lastly, we let $Y_n^*$ denote the number of consecutive
nodes on the spine, starting at the parent of $u^*$, whose non-spine children all have an odd peel number.
Assuming that $\tau(T_n, \beta_n) = \tau(T_\infty,\beta_n)$, we have the inequality
$$m(T_n) \leq \max\big( m(T_n^*), \rho(u^*)+Y_n^*, M_n' + Y_n, M_n''\big) 
\leq \max\big(m(T_n^*) + Y_n^*, 2M_n', 2Y_n, M_n''\big).\adveq$$
To explain this, we note that nodes in $S^*$ have a peel number that is at most
one more than the peel numbers of their
children on the spine. Nodes on the spine that are not in $S^*$ have a peel number that is at most one more
than the maximum peel number of any of their non-spine children (and this is bounded from above by $M_n'$).

Let $\eps>0$ be given and let $x = (1+\eps) \log n/\log\big(1/f'(1-q)\big)$. We have
\newcount\seventerms
\seventerms=\eqcount
$$\eqalign{
\pr\big\{ m(T_n)\ge x\big\} \leq \pr\big\{ &\tau(T_n,\beta_n)\neq \tau(T_\infty, \beta_n)\big\}
+ \pr\Big\{ \max_{u\in S} h(T_u) \geq \beta_n - \alpha_n\Big\}
+ \pr\Big\{Y_n^* \geq \sqrt{\log n}\Big\}\cr
&+ \pr\{Y_n\geq x/2\}
+ \pr\{ M_n' \geq x/2\}
+ \pr\{M_n'' \geq x\}
+ \pr\Big\{m(T_n^*) \geq x-\sqrt{\log n}\Big\}.\cr
}\adveq$$
As noted in our proof of the lower bound, the first two terms are $o(1)$, so we have reduced our task to
showing that the latter five terms are also $o(1)$.

Let $\zeta$ be the offspring distribution of nodes on the spine (recall that $\pr\{\zeta = i\} = ip_i$).
For a node on the spine, the probability that it is in $S^*$ is
$$\ex\big\{(1-q)^{\zeta-1}\} = \sum_{i\geq 0} ip_i (1-q)^{i-1} = f'(1-q).\adveq$$
Thus, $Y_n^*$ is a geometric random variable with parameter $1-f'(1-q)$, and hence
$$\pr\Big\{Y_n^*\geq \sqrt{\log n}\Big\} = o(1).\adveq$$
Also, $Y_n$ is bounded from above in distribution by the maximum of $\alpha_n$ independent
$\Geo\big(f'(1-q)\big)$ random variables, so that
$$\pr\{Y_n\geq x/2\} \leq \alpha_n f'(1-q)^{x/2} = o(1).\adveq$$
Next,
$$\pr\{M_n'\geq x/2\} \leq \ex\big\{|S|\big\} \pr\{R\geq x/2\}=\sigma^2\alpha_n f'(1-q)^{x/2+o(x)}=o(1)\adveq$$
and
$$\pr\{M_n''\geq x\}\leq\ex\big\{|S|\big\}\pr\{M\geq x\}=\sigma^2 \alpha_n f'(1-q)^{x/2+o(x)} = o(1).\adveq$$
This leaves us with the final term of~\refeq{\the\seventerms}. Observe that $|T_n^*| = n-\alpha_n - \sum_{u\in S}
|T_u|$, which is at most $n - \max_{u\in S} |T_u|$. Thus,
$$\eqalign{
\pr\bigg\{|T_n^*|\geq n-{n\over\log^5 n}\bigg\}&\leq\pr\bigg\{\max_{u\in S} |T_u|\leq {n\over \log^5 n}\bigg\}\cr
&= \ex\Bigg\{ \pr\bigg\{|T|\leq {n\over \log^5 n}\bigg\}^{|S|}\Bigg\} \cr
&\le \pr\bigg\{ |S|\leq {\sigma^2\alpha_n\over 2}\bigg\}
   + \bigg(1-\pr\bigg\{|T|> {n\over \log^5 n}\bigg\}\bigg)^{\sigma^2\alpha_n/2}.\cr
}\adveq$$
Since, as noted earlier, $|S|/(\sigma^2\alpha_n)\to 1$ in probability and since $\pr\big\{|T|\ge n\big\}
=\Theta(1/\sqrt n)$, we have
$$\eqalign{
\pr\bigg\{|T_n^*|\geq n-{n\over\log^5 n}\bigg\}&\leq o(1) + \exp\!\Bigg(\!\!-\!\Theta\bigg({\log^{5/2} n\over
\sqrt n}\bigg){\sigma^2\over 2}\alpha_n\Bigg)\cr
&\leq o(1) + \exp\!\bigg(\!\!-\!\Theta\Big(\sqrt{\log n}\Big)\bigg),\cr
}\adveq$$
which is $o(1)$. So
$$\pr\Big\{ m(T_n^*) \geq x-\sqrt{\log n}\Big\}\leq
\max_{1\leq k \leq n-n/\log^5 n} \pr\Big\{ m(T_n^*) \geq x - \sqrt{\log n} \Big|\, |T_n^*| = k\Big\}
+ \pr\bigg\{ |T_n^*| \geq n-{n\over \log^5 n}\bigg\}.\adveq$$
Noting that given $|T_n^*|=k$, $T_n^*$ is again a Bienaym\'e--Galton--Watson tree and letting $F_k$
be the event that there exists a node $v\in T_n$ with $|T_v|\leq k$ and $m(T_v) \geq x - \sqrt{\log n}$, we see that
$$\pr\Big\{ m(T_n^*) \geq x-\sqrt{\log n} \Big|\, |T_n^*|=k\Big\}
\leq \pr\Big\{ F_{n-n/\log^5 n}\Big\}.\adveq$$
Now define $t(v)$ to be the subtree of $v$ in the shifted preorder degree sequence $\xi_v, \xi_{v+1}, \dots, \xi_n, \xi_1, \dots \xi_{v-1}$. Let $\rho(v)$ be the peel number of the root $v$ of
$t(v)$ and let $G_{vx}$ denote the event that $\rho(v)\ge x-\sqrt{\log n}$ and $|t(v)|\le n/\log^5 n$. We have
$$\pr\{M_n\ge x\} \leq \pr\Bigl\{\bigcup_{v\in T_n} G_{vx}\Bigr\} + o(1).\adveq$$
Note that
$\max_{v\in T_n; |t(v)| \leq n/\log^5 n} \rho(v)$, is invariant under the cyclic shift of the preorder degree sequence. This
rotational invariance, by Dwass' device~\bref{11}, shows that
\newcount\peeldwass
\peeldwass=\eqcount
$$\pr\Big\{\bigcup_{v\in T_n} G_{vx}\Big\} = {\pr\big\{ \bigcup_{v\in T_n} G_{vx},\, \sum_{1\leq i\leq n}
(\xi_i - 1) = -1\big\}\over \pr\big\{ \sum_{1\leq i\leq n} (\xi_i-1) = -1\big\}},\adveq$$
where on the right-hand side, all probabilities are with respect to an i.i.d.\ sequence $\xi_1,\ldots,\xi_n$.
We bound~\refeq{\the\peeldwass} from above by
$$ \pr\Big\{\bigcup_{v\in T_n} G_{vx}\Big\} \leq
n \cdot {\pr\Big\{ \rho(1) \geq x-\sqrt{\log n},\, |t(1)|\leq n/\log^5 n,\, \sum_{i=1}^n (\xi_i-1) = -1\Big\}\over
\pr\big\{ \sum_{i=1}^n (\xi_i-1) = -1\big\}}.\adveq$$
By conditioning on the size of $t(1)$, we obtain the further bound
$$ \pr\Big\{\bigcup_{v\in T_n} G_{vx}\Big\} \leq
n\cdot \pr\Big\{\rho(1)\geq x-\sqrt{\log n}\Big\}
   \cdot {\sup_{n/\log^5 n \leq k\leq n} \pr\big\{\sum_{i=1}^k (\xi_i- 1) = 0\big\}
\over \pr\big\{\sum_{i=1}^n(\xi_i-1) = -1\big\}}.$$
By Kolchin's estimate~\bref{19}, the fraction is $\Theta(1)$, therefore,
$$ \pr\Big\{\bigcup_{v\in T_n} G_{vx}\Big\} \leq n f'(1-q)^{x+o(x)},$$
which goes to $0$ if $x = (1+\eps) \log n/ \log\bigl(1-f'(1-q)\bigr)$.\slug

If, instead of removing leaves and parents at each step, we only remove leaves, then it is clear that the
number of rounds needed to delete all nodes is simply the height of the tree. The height of random binary trees was
studied by P.~Flajolet and A.~Odlyzko, who showed that in this case, $H_n/\sqrt n$
converges in law to a theta distribution~\bref{12}.
Earlier, it was shown by N.~G.~de~Bruijn, D.~E.~Knuth, and S.~O.~Rice
that the expected height of a random planted plane tree is $\sqrt{\pi n}+O(1)$. It is interesting that
deleting only leaves from $T_n$ at each step requires $\Theta\bigl(\sqrt n\bigr)$
rounds of deletion, but deleting leaves
and their parents causes the number of rounds to decrease to $\Theta(\log n)$.

\medskip
\boldlabel Examples. We apply Theorem~{\the\peelthm} to calculate explicit asymptotics of the maximum peel number for the various families of trees mentioned earlier.
\medskip
\item{i)} {\it Flajolet's $t$-ary trees:} We have $f'(1-q) = 1-q$ and thus $M_n/\log n \to 1/\log(1/(1-q))$ in probability. As $t$ gets large, $q$ approaches $1-1/t$, so that the limit of $M_n/\log n$ is approximately
$1/\log t$ for large $t$.
For the case of full binary trees when $t=2$, recall that $q = 2 - \sqrt{2}$ and thus
$M_n/\log n \to -1/\log(\sqrt{2} - 1)$ in probability.
\smallskip
\item{ii)} {\it Cayley trees:} In this case, $f'(1-q) = e^{-q}$ and hence $M_n/\log n \to 1/q$ in probability;
we know from earlier that $q=W(1)$, so $1/q\approx 1.763223$.
\smallskip
\item{iii)} {\it Planted plane trees:} We calculate $f'(1-q) = 1/(q+1)^2$. Recalling that $q = 1/\varphi$ where $\varphi = (\sqrt{5}-1)/2$ is the golden ratio, we have in probability $M_n/\log n \to 1/\varphi^2 \approx 0.381966$.
\smallskip
\item{iv)} {\it Motzkin trees:} The derivative $f'(1-q) = (3-2q)/3$ and substituting $q = 3 - \sqrt{6}$ we get $M_n/\log n \to 1/(\log 3 - \log (2\sqrt{6}-3)) \approx 2.186769$.
\smallskip
\item{v)} {\it Binomial trees:} In this case, we have $f'(1-q) = (1-q/d)^{d-1}$.
As $d \to \infty$, it is clear to see that $M_n/\log n \to 1/W(1)$ in probability, matching the earlier calculation for Cayley trees above. For the special case $d = 2$ of Catalan trees, we have $f'(1-q) = 1-q/2$ and thus $M_n/\log n \to -1/\log(\sqrt{3}-1)$. This constant is greater than that we obtain for full binary trees above, which is consistent with intuitive reasoning about the maximal peel numbers of these trees.
\medskip\goodbreak

\advsect Distribution of the leaf-height

We now repeat the treatment given in Section~{\the\peelnumsect}, but this time
for the distribution $(\ell_i)_{i\geq 0}$, where for each $i\geq 0$,
$\ell_i$ denotes the probability that the root of an unconditional
Bienaym\'e--Galton--Watson tree has leaf-height equal to
$i$.
Observe that $\ell_0$ is exactly the probability $p_0$ that the root has no children and in general,
for a node $u$ with children $\Gamma_u$ the leaf-height $\leaf(u)$ is
$$\leaf(u) = \cases{0,&if $u$ is a leaf;\cr \min_{v\in \Gamma_u}\leaf(v) + 1,& otherwise.}\adveq$$
We define $\lp_i$ and $\lm_i$ analogously to $\rplus_i$ and $\rminus_i$:
$$\lp_i = \sum_{j\geq i}\ell_j\qquad\hbox{and}\qquad \lm_i = \sum_{j=0}^i \ell_j\,;\adveq$$
since $(\ell_i)_{i\geq 0}$ defines a distribution, $\lp_{i+1} + \lm_i = 1$ for every $i\geq 0$.
Letting $E_i$ be the event that all the children of the root have leaf-height
at least $i$, we have, for $i\geq 1$, $\ell_i = \pr\{E_{i-1}\} - \pr\{E_i\}$. We can then compute
$$\ell_1 = 1-\ex\big\{(1-\ell_0)^\xi\big\} = 1-f(1-\ell_0) = 1-f(1-p_0).\adveq$$
and, in general, for $i\geq 1$,
$$\eqalign{
\ell_{i+1} &= \pr\{E_i\} - \pr\{E_{i+1}\} \cr
&= \ex\big\{(\lp_i)^\xi\big\} - \ex\big\{(\lp_{i+1})^\xi\big\}\cr
&= f(\lp_i) - f(\lp_{i+1})\cr
&= f(1-\lm_{i-1}) - f(1-\lm_i).\cr
}\adveq$$
By convexity of $f$, we see that $\ell_i$ is nonincreasing, and this formula provides a fast method to compute
the $\ell_i$ recursively.
The following lemma describes the behaviour of $\ell_i$ as $i$ gets large.

\proclaim Lemma \advthm. Let $T$ be a Bienaym\'e--Galton--Watson tree with offspring distribution
$\xi\sim(p_i)_{i\geq 0}$. If $p_1 \neq 0$, then $\ell_i = (p_1 + o(1))^i$.
Otherwise if $p_1 = 0$ and $\kappa = \min\{i>1:p_i\neq 0\}$, then
$$\log \ell_i = \Theta(\kappa^i)\adveq$$
as $i\to \infty$.
\goodbreak

\proof The recursive formula above is our starting point.
Expanding $f$ as a power series, for $i\geq 1$ we have, by our choice of $\kappa$,
$$\eqalign{
\ell_{i+1} &= \sum_{j=0}^\infty p_j\big((\lp_i)^j - (\lp_{i+1})^j\big) \cr
&= 0 + p_1(\lp_i - \lp_{i+1}) + \sum_{j\geq\kappa}p_j\big((\lp_i)^j - (\lp_{i+1})^j\big)\cr
&= p_1 \ell_i + p_\kappa \ell_i\big((\lp_i)^{\kappa-1}
  + (\lp_i)^{\kappa - 2}(\lp_{i+1})^{1} + \dots + (\lp_{i+1})^{\kappa - 1}\big)
   + \sum_{j>\kappa}p_j\big((\lp_i)^j - (\lp_{i+1})^j\big).\cr
&\leq p_1\ell_i + \kappa p_\kappa \ell_i(\lp_i)^{\kappa-1} + \sum_{j>\kappa} jp_j\ell_i(\lp_i)^{j-1}\cr
&\leq p_1\ell_i + \ell_i(\lp_1)^{\kappa-1}\Big(\sum_{j\geq \kappa} jp_j\Big)\cr
&=p_1\ell_i + \ell_i(1-p_0)^{\kappa-1}(1-p_1).\cr
}\adveq$$
Letting $\alpha = p_1 + (1-p_1)(1-p_0)^{\kappa-1}<1$, we have $\ell_{i+1} \leq \ell_i\alpha$. Hence $\ell_{i+1}\leq
\ell_1\alpha^i$ and therefore $\ell_i\to 0$ as $i\to\infty$.

Let $\eps>0$ and pick $n_\eps$ large enough such that $\lp_i\leq \eps$ for all $i\geq n_\eps$. When $p_1\neq
0$, we have $\ell_ip_1\leq \ell_{i+1} \leq \ell_i(p_1+\eps^{\kappa-1})$, so we immediately conclude that $\ell_i =
\big(p_1+o(1)\big)^i$. If $p_1 = 0$, then
$$\eqalign{
\ell_{i+1}&\leq\kappa p_\kappa \ell_i(\lp_i)^{\kappa-1}
  + \sum_{j>\kappa}jp_j\ell_i(\lp_i)^{\kappa-1}(\lp_i)^{j-\kappa}\cr
&\leq \kappa p_\kappa \ell_i (\lp_i)^{\kappa-1} \Big(1+\sum_{j>\kappa} jp_j\eps^{j-\kappa}\Big)\cr
&\leq \kappa p_\kappa \ell_i (\lp_i)^{\kappa-1} \bigg(1+{\eps\over 1-\eps}\bigg)\cr
&\leq \kappa p_\kappa \ell_i (\ell_i)^{\kappa-1}
      \Big(\sum_{j=0}^\infty \alpha^j\Big)^{\kappa -1}\bigg({1\over1-\eps}\bigg)\cr
&= {\kappa p_\kappa\over (1-\alpha)^{\kappa - 1}(1-\eps)}{\ell_i}^\kappa.\cr
}\adveq$$
From this and the fact that $\ell_i\to 0$, we see that for some positive constants $c_1$, $c_2<1$, and $c_3$,
$$\l_i\leq c_1{c_2}^{\kappa^{i-c_3}}\adveq$$
for all $i\geq c_3$. We also have
$$\l_{i+1} \geq \kappa p_\kappa \l_i(\lp_i)^{\kappa-1}\geq \kappa p_\kappa{\l_i}^k,\adveq$$
so that for some positive constants $c_1'$, $c_2'<1$, and $c_3'$,
$$\l_i\geq c_1'{c_2'}^{\kappa^{i-c_3'}}\adveq$$
for all $i\geq c_3'$. This proves that $\log \l_i = \Theta(\kappa^i)$. We finish the proof by noting
that $\lp_i$ can be bounded in a similar manner.\slug

\advsect Asymptotics of the leaf-height

In this section, we will describe the asymptotic behaviour of
the leaf-height of a tree $T_n$ (recall that this is the maximum of $\leaf(v)$, taken over all the nodes
$v\in T_n$). The result depends on whether $p_1$ is zero or nonzero, and we have
split this into two lemmas, since both of the proofs are rather involved.

\newcount\ponenonzero
\ponenonzero=\thmcount
\proclaim Lemma \advthm. Let $L_n$ be the leaf-height of $T_n$,
a conditional Bienaym\'e--Galton--Watson tree on $n$ nodes with
offspring distribution $\xi\sim (p_i)_{i\geq 0}$. If $p_1\neq 0$, then
$${L_n\over \log n} \to {1\over \log(1/p_1)}\adveq$$
in probability.

\proof Let $Y_n$ be the length of the longest string, oriented away from the root, of nodes
of degree one in $T_n$. Clearly $L_n\geq Y_n$, so we will first show that for $\eps>0$,
$$\pr\big\{Y_n<(1-\eps) \log n / \log(1/p_1)\big\}\to 0.$$
Now, $\pr\{Y_n<x\}$ is the probability that a string of 1s appears
in the preorder degree sequence of the tree $(\xi_1, \xi_2, \ldots,\xi_n)$, given that the sequence is of
length $n$ and that the sequence does, in fact, define a tree; as we have used previously, this latter probability is $\Theta(n^{-3/2})$,
from Dwass~\bref{11}. So letting $Y_n(\xi_1,\xi_2,\ldots,\xi_n)$ be the length of the longest subsequence of
1s in the preorder degree sequence, we have
$$\pr\{Y_n<x\} = \Theta(n^{-3/2})\pr\big\{Y_n(\xi_1,\xi_2,\ldots,\xi_n)<x\big\}.\adveq$$
We divide the sequence into $n/x$ subsequences of length $x$ each and let $E_i$ be the event that the
$i$th subsequence {\it does not} consist only of 1s. Then
$$\pr\big\{Y_n(\xi_1,\xi_2,\ldots,\xi_n)<x\big\} = \pr\bigg\{\bigcup_{i=1}^{n/2} E_i\bigg\}
  =\pr\{E_i\}^{n/2}.\adveq$$
Since $\pr\{E_i\} = 1-{p_1}^x$ for all $i$,
$$\pr\big\{Y_n(\xi_1,\xi_2,\ldots,\xi_n)<x\big\}
= (1-{p_1}^x)^{n/x} \leq \exp\!\bigg(\!\!-\!{n{p_1}^x\over x}\bigg).\adveq$$
When $x = (1-\eps) \log n / \log(1/p_1)$, this is equal to $\exp(-n^\eps/x)$, so
$$\pr\big\{L_n < (1-\eps) \log n / \log(1/p_1)\big\} \leq \pr\big\{Y_n < (1-\eps) \log n / \log(1/p_1)\big\}
\leq \Theta(n^{3/2})e^{-\Theta(n^\eps/\log n)},\adveq$$
which goes to $0$ as $n\to\infty$.

To tackle the upper bound, it will be helpful for us to reorder the degrees into level (also called
breadth-first) ordering and to consider the following random variable. Arrange the level-order degree sequence
$\xi_1,\xi_2,\ldots,\xi_n$ in a cycle, and let
$Z_n$ be the longest string of consecutive non-zero numbers in this cyclic ordering. Clearly the probability
that no sub-cycle of length $x$ of this ordering consists only of zeroes is at most $n(1-p_0)^x$.
So letting $A$ be the event that $(\xi_1,\xi_2,\ldots,\xi_n)$
defines a tree, we can crudely bound $\pr\{Z_n\geq x\}$ by
$$\pr\{Z_n\geq x\} = {\pr\{Z_n\geq x,\,A\}\over \pr\{A\}} \leq
\Theta(n^{3/2})n(1-p_0)^x=\Theta(n^{5/2})(1-p_0)^x,\adveq$$
which goes to zero if $c > (5/2)/\log\big(1/(1-p_0)\big)$ and $x$ is set to $c\log n$. By symmetry, if $Z_n$ is the
longest string of nonzeroes in the {\it preorder} listing, then the same result holds, that is,
$$\pr\bigg\{Z_n\geq {3\log n\over \log\big(1/(1-p_0)\big)}\bigg\}\to 0.\adveq$$
Note that $L_n\leq Z_n$. Now for $1\leq \Delta\leq n$ let $L_n(\xi_1,\xi_2,\ldots,\xi_\Delta)$ be the
smallest leaf depth if we start constructing a tree using degrees $\xi_1,\xi_2,\ldots,\xi_\Delta$, in preorder.
Two situations can occur: either $\xi_1,\xi_2,\ldots,\xi_\Delta$ defines at least one tree in a possible
forest, or $\xi_1,\xi_2,\ldots,\xi_\Delta$ defines an incomplete tree. In the former case,
$L_n(\xi_1,\xi_2,\ldots,\xi_\Delta)$ is the leaf-height of the first completed tree; in the latter, set
$L_n(\xi_1,\xi_2,\ldots,\xi_\Delta) = 0$. Note that if $Z_n\leq \Delta$, then
$L_n(\xi_1,\xi_2,\ldots,\xi_\Delta)\leq \Delta$. For a sequence $\xi_1,\xi_2,\ldots,\xi_n$ of degrees, we define
$$L_{ni} = L_n(\xi_i,\xi_{i+1}\ldots,\xi_{i+\Delta-i}),\adveq$$
where addition in the indices is taken modulo $n$. If $Z_n\leq \Delta$, note that $L_n\leq \max_{1\leq i\leq
n} L_{ni}$. So we have
$$\eqalign{
\pr\{L_n>x\} &\leq \pr\{L_n>x,\,Z_n\leq \Delta\} + \pr\{Z_n>\Delta\} \cr
&\leq \pr\Big\{\max_{1\leq i\leq n} L_{ni}>x,\,Z_n\leq \Delta\Big\} + \pr\{Z_n>\Delta\}\cr
&\leq \pr\Big\{\max_{1\leq i\leq n} L_{ni}>x\Big\} + \pr\{Z_n>\Delta\}.\cr
}\adveq$$
The second term is $o(1)$ if we pick $c > (5/2)/\log\big(1/(1-p_0)\big)$ as before and set $\Delta = c\log n$.
In the first term, the maximum is invariant under rotations of $(\xi_1,\xi_2,\ldots,\xi_n)$, so we use may use a
version of the cycle lemma~\bref{11}, obtaining
$$\eqalign{
\pr\Big\{\max_{1\leq i\leq n} L_{ni}>x\Big\} &=
 {\pr\big\{\max_{1\leq i\leq n} L_{ni}>x,\,\sum_{i=1}^n\xi_i=n-1\big\}\over\pr\big\{\sum_{i=1}^n\xi_i=n-1\big\}}\cr
&\leq
 {\pr\big\{\max_{1\leq i\leq n} L_{ni}>x,\,\sum_{i=\Delta+1}^n\xi_i=n-1-\sum_{i=1}^\Delta\xi_i\big\}\over
 \Theta(n^{-1/2})}\cr
&\leq O(n^{3/2})\sup_\ell\pr\bigg\{L_{n1}>x,\,\sum_{i=\Delta+1}^n\xi_i=\ell\bigg\}\cr
&= O(n^{3/2})\cdot\pr\{L_{n1}>x\}\cdot\sup_\ell\,\pr\bigg\{\sum_{i=\Delta+1}^n\xi_i=\ell\bigg\}.\cr
}\adveq$$
Rogozin's inequality~\bref{24} tells us that
$$\sup_\ell\,\pr\Biggl\{\sum_{i=\Delta+1}^n\xi_i=\ell\Biggr\}
  \leq {\gamma\over\sqrt{1-\Pi}}\cdot{1\over\sqrt{n-\Delta}},\adveq$$
where $\gamma$ is a universal constant and $\Pi = \sup_j p_j$. So if $L(T)$ is the leaf-height of the root of
an unconditional Bienaym\'e--Galton--Watson tree, then
$$\pr\{L_n>x\} \leq O(n) \pr\{L_{n1}>x\}\leq O(n)\pr\{L(T)>x\}\leq O(n)\lp_x\leq O(n)\big(p_1+o(1)\big)^x,\adveq$$
which goes to 0 when $x = (1+\eps)\log n/\log(1/p_1)$.\slug

The next lemma handles the other case, in which $p_1$ is zero.

\newcount\ponezero
\ponezero=\thmcount
\proclaim Lemma \advthm. Let $L_n$ be the leaf-height of $T_n$, a conditional
Bienaym\'e--Galton--Watson tree on $n$ nodes with
offspring distribution $\xi\sim (p_i)_{i\geq 0}$. Let $\kappa = \min\{i>1: p_i\neq 0\}$. If $p_1 = 0$, then
$${L_n\over \log\log n} \to {1\over \log\kappa}\adveq$$
in probability.

\proof Let $\eps>0$, $A$ be the event that
$(\xi_1,\xi_2,\ldots,\xi_n)$ forms a tree and let $L_n(\xi_1,\xi_2,\ldots,\xi_n)$ be as in the previous lemma.
If $L(T)$ is the leaf-height of the root of an unconditional Bienaym\'e--Galton--Watson tree $T$, we have,
by the cycle lemma,
$$\eqalign{
\pr\{L_n\geq x\} &= {\pr\big\{L_n(\xi_1,\xi_2,\ldots,\xi_n)\geq x,\, A\big\}\over \pr\{A\}} \cr
&\leq \Theta(n^{3/2}) \pr\big\{L_n(\xi_1,\xi_2,\ldots,\xi_n)\geq x\big\} \cr
&\leq \Theta(n^{3/2}) n\pr\{L(T)\geq x\} \cr
&\leq \Theta(n^{5/2}) c_1 c_2 ^{\kappa^{x-c_3}},\cr
}\adveq$$
for some positive constants $c_1$, $c_2<1$, and $c_3$. When $x = (1+\eps)\log_\kappa\log n$, this
is $\Theta(n^{5/2}) c_2^{\Theta((\log n)^{1+\eps})}$, which goes to $0$ as $n\to\infty$.

Let $T_\infty$ be Kesten's infinite tree. Our proof of the lower bound uses the fact that
$$\pr\big\{\tau(T_\infty, n^{1/3})\neq \tau(T_n,n^{1/3})\big\}\to 0\adveq$$
as $n\to \infty$. Let $U$ be the set of all unconditional Bienaym\'e--Galton--Watson
trees $T$ rooted less than $n^{1/4}$ of
the way down the spine. Let $h(T)$ denote the height of an unconditional Bienaym\'e--Galton--Watson tree; the
probability that one of these trees has height greater than $n^{1/3}/2$ is bounded above by
$$\ex\big\{|U|\big\}\pr\big\{h(T)>n^{1/3}\big\} \leq \sigma^2n^{1/4} \bigg({2+o(1)\over
\sigma^2n^{1/3}/2}\bigg)\sim {4\over n^{1/12}},\adveq$$
which goes to $0$.
Here we used the fact that $\ex\{\zeta\} = \sigma^2+1$ and applied Kolmogorov's result
(see~\bref{20} and~\bref{2}) about
the height of a Bienaym\'e--Galton--Watson tree.
If all the heights above are at most $n^{1/3}/2$, then $T_n^*$, the tree
obtained by taking the spine up to level $n^{1/4}$ and all hanging unconditional trees up to that point, is
a subtree of $\tau(T_\infty, n^{1/3})$, since $n^{1/4} + n^{1/3}/2 < n^{1/3}$.
\midinsert
\vskip5pt
$$\epsfbox{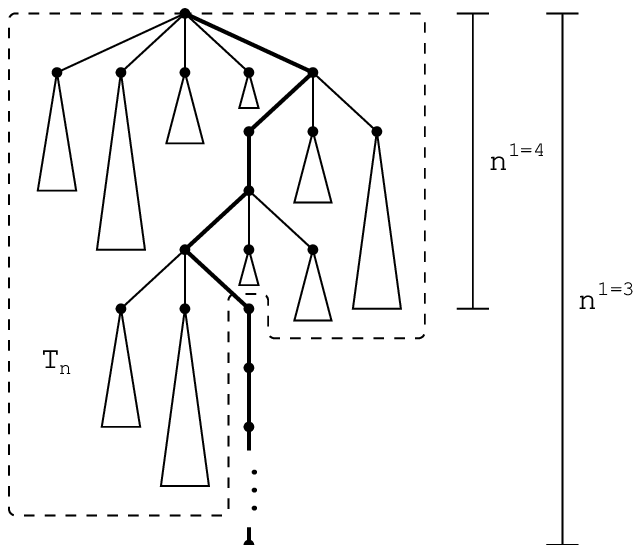}$$
\caption{The construction in the lower bound. Note that none of the unconditional
trees in $T_n^*$ reach level $n^{1/3}$.}
\endinsert
For every tree $t\in U$,
let $E_t$ be the event that the leaf-height of the root of $t$ is less than $x$ and let $E_T$ be the same
event for an unconditional Bienaym\'e--Galton--Watson tree (since each $t$ in $U$ is such a tree,
there is no real moral
distinction between these events). We have
$$\eqalign{
\pr\{L_n<x\}&\leq\pr\big\{\tau(T_\infty)\neq (T_n,n^{1/3})\big\} + \pr\big\{T_n^*\neq (T_n,n^{1/3})\big\}
  + \pr\Big\{\bigcap_{t\in U} E_t\Big\} \cr
&\leq o(1) + \ex\big\{ \pr\{E_T\}^{|U|}\big\}\cr
&\leq o(1) + \ex\big\{ \big(1-c_1'\cdot c_2'^{\kappa^{x-c_3'}}\big)^{|U|}\big\}\cr
&\leq o(1) + \ex\big\{ \big(1-c_1'\cdot c_2'^{\kappa^{x-c_3'}}\big)^{n^{1/4}\sigma^2/2}\big\}
  +\pr\bigg\{|U|<{n^{1/4}\sigma^2\over 2}\bigg\},\cr
}\adveq$$
for some $c_1',c_2',c_3'$ positive, $c_2'<1$. Take $x = (1-\eps)\log_\kappa\log n$. Letting
$\zeta_1,\zeta_2,\ldots,\zeta_n$ be independent and distributed as $\zeta$, we find that
$$\pr\{L_n<x\}\leq o(1) + \exp\!\big(-\Theta(n^{1/4}e^{-\Theta((\log n)^{1-\eps})})\big)
  + \pr\bigg\{\sum_{i=1}^{n^{1/4}} (\zeta_i - 1) < {n^{1/4}\sigma^2\over 2}\bigg\}.\adveq$$
Since $\ex\{\zeta-1\} = \sigma^2$, by the law of large numbers, this entire expression is $o(1)$.\slug

It is important to note that the node with maximum leaf-height in a tree usually is not the root. We have
the following result for the distribution of the leaf-height of the root of a conditional
Bienaym\'e--Galton--Watson tree $T_n$.

\proclaim Lemma \advthm. Let $L_n'$ denote the leaf-height of the root of $T_n$,
a conditional Bienaym\'e--Galton--Watson tree of size
$n$ and offspring distribution $\xi$.
Let $f$ be the generating function of this distribution. Then the probability distribution of $L_n'$ is given by
$$\lim_{n\to\infty}\pr\{L_n' = i\} = \prod_{j=0}^{i-1}f'(\lp_j).\adveq$$

\proof For a node $v$ on the spine of Kesten's infinite tree $T_\infty$, let $H^*$ be the
leaf-height of the spine-child of $v$, and $H(1),H(2),\ldots,H(\zeta-1)$ be the leaf-heights of the $\zeta-1$
independent unconditional
Bienaym\'e--Galton--Watson trees spawned by $v$. Then the leaf-height of a node $v$ on the spine is
$$1+\min_{w\in\Gamma_v} \leaf(w) = 1+\min \big(H(1), H(2),\ldots,H(\zeta-1)\big).\adveq$$
This defines a Markov chain on the positive integers that proceeds up the spine. The state $H^*$ (which is
just a positive integer indicating the leaf-height of the node on the spine), is taken to the state
$$1+\min\!\big(H^*, H(1), H(2), \ldots,H(\zeta-1)\big)\adveq$$
in one step of the Markov chain; here all $H_i$ have distribution $(\l_i)_{i\geq 0}$. Let $H^{**}$ be the
limit stationary random variable of this Markov chain. That the limit exists and that the chain is positive
recurrent follows from the fact that at each step, there is a positive probability that the next state is $1$.
This happens when $\zeta>1$ and one of $H(i)$ is 0. In fact, $H^{**}$ is the unique solution of the
distributional identity
$$H^{**}\;{\buildrel{\cal L}\over=}\;1 + \min\!\big(H^{**}, H(1), H(2),\ldots,H(\zeta-1)).\adveq$$
Broutin, Devroye, and Fraiman showed that under a coalescence condition (satisfied here), the limit of the root
value of $T_n$ tends in distribution to the stationary random variable for Kesten's spinal Markov
chain~\bref{7}. Thus
$L_n'\to H^{**}$ in distribution.

We can describe the distribution of $H^{**}$ more explicitly. For convenience, let $\lss_i = \pr\{H^{**}=i\}$.
For $i\geq 1$,
$$\lss_i = \pr\{H_j \geq i+1\ \hbox{for all}\ 1\leq j\leq \zeta -1\}\pr\{H^{**}\geq i-1\}.\adveq$$
This means that
$$ \lss_i = \lss_{i-1}\ex\big\{(\lp_{i-1})^{\zeta-1}\big\}=\lss_{i-1}\sum_{j\geq 1}jp_j (\lp_{i-1})^{j-1},\adveq$$
and we can rewrite this in terms of the generating function $f(s)$ of $\xi$ as
$$\lss_i = \lss_{i-1} f'(\lp_{i-1}) = \lss_0\prod_{j=0}^{i-1}f'(\lp_j) = \prod_{j=0}^{i-1}f'(\lp_j),\adveq$$
proving the lemma.\slug

Lastly, we can obtain a random variable by taking the leaf-height of a node chosen uniformly at random from
$T_n$. Its distribution is asymptotically the same as the leaf-height of the root of an unconditional
Bienaym\'e--Galton--Watson tree.

\proclaim Lemma \advthm. Let $L_n''$ be a random variable obtained by taking the leaf-height of a node chosen
uniformly at random from a conditional Bienaym\'e--Galton--Watson tree $T_n$. We have
$$\lim_{n\to\infty}\pr\{L_n'' = i\} = \l_i\adveq$$
for all $i\geq 0$.

\proof By Aldous's theorem~\bref{1}, if $T_n^*$ denotes the subtree of $T_n$ rooted at a uniformly selected
random node, then for all trees $t$,
$$\lim_{n\to\infty} \pr\{T_n^*=t\} = \pr\{T = t\},\adveq$$
where $T$ is the unconditional Bienaym\'e--Galton--Watson tree. The result is immediate.\slug

\boldlabel Examples. We now apply Lemmas~{\the\ponenonzero} and~{\the\ponezero} to compute the maximum leaf-height
(asymptotically in probability)
for common families of trees.
These results, along with the independence numbers and maximum peel numbers we computed earlier, are
collected in Table 1. In the table, $W$ denotes the Lambert function and $\varphi = (\sqrt{5}-1)/2$ is the golden ratio.
\midinsert
$$\vcenter{\vbox{
\vskip -15pt
\centerline{\smallheader Table 1}
\smallskip
\centerline{\eightpoint
ASYMPTOTIC VALUES OF PARAMETERS FOR CERTAIN FAMILIES OF TREES}
}}$$
$$\vcenter{\vbox{
\eightpoint
\tabskip=.5em plus2em minus.6em
    \halign{
        \hfil$\displaystyle{#}$\hfil & \hfil$\displaystyle{#}$\hfil & \hfil$\displaystyle{#}$\hfil &
        \hfil$\displaystyle{#}$\hfil & \hfil$\displaystyle{#}$\hfil \cr
        \noalign{\hrule}
        \noalign{\medskip}
        \hbox{Family} & I_n & M_n & L_n \cr
        \noalign{\medskip}
        \noalign{\hrule}
        \noalign{\medskip}
        {\hbox{Full binary} \atop (\Uni\{0,2\})}  & (2 - \sqrt 2)n 
             &  {\log n\over \log(1/(\sqrt 2-1))} &  \log_2\log n \cr
        \noalign{\medskip}
        {\hbox{Flajolet $t$-ary}\atop (p_0 = 1-1/t;\;p_t= 1/t)}  & \Bigl(1-{1+o_{t\to\infty}(1)\over t}\Bigr)n
             & \sim_{t\to\infty} \log_t n & \log_t \log n  \cr
        \noalign{\medskip}
        {\hbox {Cayley} \atop (\Pos(1))} & W(1)n & {\log n/ W(1)} & \log n\cr
        \noalign{\medskip}
        {\hbox{Planted plane} \atop (\Geo(1/2))} & {n\over \varphi} & {\log n\over\varphi^2} & \log_4 n \cr
        \noalign{\medskip}
        {\hbox{Motzkin} \atop (\Uni\{0,1,2\})} & (3-\sqrt 6)n 
             & {\log n\over \log 3 - \log (2\sqrt 6-3)} & \log_3 n\cr
        \noalign{\medskip}
        {\hbox{Catalan} \atop (\Bin(2,1/2))} & (4-2\sqrt 3)n & {\log n\over \log (1/(\sqrt 3-1))} & \log_2 n
              \cr
        \noalign{\medskip}
        {\hbox{Binomial} \atop (\Bin(d,1/d))} & \sim_{d\to\infty} W(1)n
             & \sim_{d\to\infty} {\log n\over W(1)} & {\log n\over (d-1)\log(1/(1-1/d))}\cr
        \noalign{\medskip}
        \noalign{\hrule}
    }
}}$$
\endinsert

\medskip
\item{i)} {\it Flajolet's $t$-ary trees:} For $t\geq 2$, we have $p_1 = 0$ here and $\kappa = t$, so we have
$L_n/\log\log n \to 1/\log t$ in probability, by Lemma~{\the\ponezero}.
\smallskip
\item{ii)} {\it Cayley trees:} In this case, $p_1 = 1/e$, so Lemma~{\the\ponenonzero} gives us
$L_n/\log n \to 1$ in probability.
\smallskip
\item{iii)} {\it Planted plane trees:} For these trees, $p_1 = 1/4$, so we have
$L_n/\log n \to 1/\log 4$ in probability, by Lemma~{\the\ponenonzero}.
\smallskip
\item{iv)} {\it Motzkin trees:} This family has $p_1 = 1/3$ and $L_n/\log_n \to 1/\log 3$ in probability.
\smallskip
\item{v)} {\it Binomial trees:} For a parameter $d\geq 2$, we have $p_1 = (1-1/d)^{d-1}$, which means that
$\L_n/\log n \to 1/\bigl((1-d)\log(1-1/d)\bigr)$ in probability.
As $d\to\infty$, the denominator approaches $1$, which
gives the same leaf-height as the case of Cayley trees. In the special case when $d=2$, we have the
Catalan trees, for which $p_1 = 1/2$ and $L_n/\log n \to 1/\log 2$ in probability.
\medskip\goodbreak

\advsect Further directions

The definition of the peel number and our characterization of its asymptotic growth fully describes the running
time of Algorithm~I, mentioned in the introduction, which computes the layered independent set.
It would be interesting to consider the runtime of the more general Algorithm~P,
described in Section~{\the\spathsect}. To this end, we define
{\it higher-order peel numbers} as follows. Algorithm P generates an $(r+1)$-path vertex cover by repeatedly
deleting subtrees with height exactly $r$ (and marking their roots). If a node $u$ is deleted in the $m$th
iteration of the loop and is at depth $i$ of the subtree that is deleted, then its {\it peel number of
order $r$} (or {\it $r$th order peel number}) is $mk-i$.
Note that the loop counter $m$ should start at $1$ and we have $0\leq i\leq r$. By this
definition, the peel number we studied in this paper is simply the first order peel number. To determine
the runtime of Algorithm~P, one should in principle be able to approach the higher order peel numbers in the
same way we approached the case $r=1$ in Sections~{\the\peelnumsect} and~{\the\peelnumasymp}. However,
even in this case one had to handle the even and odd cases separately, and we anticipate that the analysis
of higher-order peel numbers will require careful reasoning with respect to congruence modulo $r+1$.

\section Acknowledgements

All three authors are supported by the Natural Sciences and Engineering Research Council of Canada (NSERC).
We also thank Anna Brandenberger for wonderful discussions and for her valuable technical feedback.
Last but not least, we owe a debt of gratitude to the anonymous referees for their careful reading of
our paper and for their detailed and insightful comments.

\def\l{\char32l}
\def\ke{\eob}

\section References

\parskip=0pt
\hyphenpenalty=-1000 \pretolerance=-1 \tolerance=1000
\doublehyphendemerits=-100000 \finalhyphendemerits=-100000
\frenchspacing
\def\bref#1{[#1]}
\def\beginref{\noindent}
\def\endref{\medskip}
\vskip\parskip

\beginref
\parindent=20pt\item{\bref{1}}
\hldest{xyz}{}{bib1}%
David Aldous,
``Asymptotic fringe distributions for general families of random trees,''
{\sl Annals of Applied Probability}\/
{\bf 1}
(1991),
228--266.
\endref
\beginref
\parindent=20pt\item{\bref{2}}
\hldest{xyz}{}{bib2}%
Gerold Alsmeyer,
{\sl Galton--Watson Processes,} Course notes at the University of M\"unster, 2008.
\endref
\beginref
\parindent=20pt\item{\bref{3}}
\hldest{xyz}{}{bib3}%
Krishna Athreya
and Peter Ney,
{\sl Branching Processes}
(Berlin:
Springer Verlag,
1972).
\endref
\beginref
\parindent=20pt\item{\bref{4}}
\hldest{xyz}{}{bib4}%
Cyril Banderier,
Markus Kuba,
and Alois Panholzer,
``Analysis of three graph parameters for random trees,''
{\sl Random Structures and Algorithms}\/
{\bf 35}
(2009),
42--69.
\endref
\beginref
\parindent=20pt\item{\bref{5}}
\hldest{xyz}{}{bib5}%
Mikl\'os B\'ona,
``$k$-protected vertices in binary search trees,''
{\sl Advances in Applied Mathematics}\/
{\bf 53}
(2014),
1--11.
\endref
\beginref
\parindent=20pt\item{\bref{6}}
\hldest{xyz}{}{bib6}%
Bo\v{s}tjan Bre\v{s}ar,
Marko Jakovac,
J\'an Katreni\v{c},
Gabriel Semani\v{s}in,
and Andrej Taranenko,
``On the vertex $k$-path cover,''
{\sl Discrete Applied Mathematics}\/
{\bf 161}
(2013),
1943--1949.
\endref
\beginref
\parindent=20pt\item{\bref{7}}
\hldest{xyz}{}{bib7}%
Nicolas Broutin,
Luc Devroye,
and Nicolas Fraiman,
``Recursive functions on conditional Galton--Watson trees,''
{\sl Random Structures and Algorithms}\/
{\bf 57}
(2020),
304--316.
\endref
\beginref
\parindent=20pt\item{\bref{8}}
\hldest{xyz}{}{bib8}%
Gi-Sang Cheon
and Louis Welles Shapiro,
``Protected points in ordered trees,''
{\sl Applied Mathematics Letters}\/
{\bf 21}
(2008),
516--520.
\endref
\beginref
\parindent=20pt\item{\bref{9}}
\hldest{xyz}{}{bib9}%
Keith Copenhaver,
``$k$-protected vertices in unlabeled rooted plane trees,''
{\sl Graphs and Combinatorics}\/
{\bf 33}
(2017),
347--355.
\endref
\beginref
\parindent=20pt\item{\bref{10}}
\hldest{xyz}{}{bib10}%
Rosena Ruo Xia Du
and Helmut Prodinger,
``On protected nodes in digital search trees,''
{\sl Applied Mathematics Letters}\/
{\bf 25}
(2012),
1025--1028.
\endref
\beginref
\parindent=20pt\item{\bref{11}}
\hldest{xyz}{}{bib11}%
Meyer Dwass,
``The total progeny in a branching process,''
{\sl Journal of Applied Probability}\/
{\bf 6}
(1969),
682--686.
\endref
\beginref
\parindent=20pt\item{\bref{12}}
\hldest{xyz}{}{bib12}%
Philippe Flajolet
and Andrew Odlyzko,
``The average height of binary trees and other simple trees,''
{\sl Journal of Computer and System Sciences}\/
{\bf 25}
(1982),
171--213.
\endref
\beginref
\parindent=20pt\item{\bref{13}}
\hldest{xyz}{}{bib13}%
Michael Fuchs,
Cecilia Holmgren,
Dieter Mitsche,
and Ralph Neininger,
``A note on the independence number, domination number and related parameters of random binary search trees and random recursive trees,''
{\sl Discrete Applied Mathematics}\/
{\bf 292}
(2021),
64--71.
\endref
\beginref
\parindent=20pt\item{\bref{14}}
\hldest{xyz}{}{bib14}%
Bernhard Gittenberger,
Zbigniew Go\l\ke biewski,
Isabella Larcher,
and Ma\l gorzata Sulkowska,
``Protection numbers in simply generated trees and P\'olya trees,''
{\sl arXiv preprint 1904.03519}\/
(2019).
\endref
\beginref
\parindent=20pt\item{\bref{15}}
\hldest{xyz}{}{bib15}%
Svante Janson,
``Simply generated trees, conditioned Galton--Watson trees, random allocations and condensation,''
{\sl Probability Surveys}\/
{\bf 9}
(2012),
103--252.
\endref
\beginref
\parindent=20pt\item{\bref{16}}
\hldest{xyz}{}{bib16}%
Svante Janson,
``Asymptotic normality of fringe subtrees and additive functionals in conditioned Galton--Watson trees,''
{\sl Random Structures and Algorithms}\/
{\bf 48}
(2016),
57--101.
\endref
\beginref
\parindent=20pt\item{\bref{17}}
\hldest{xyz}{}{bib17}%
G\"otz Kersting,
``On the height profile of a conditioned Galton--Watson tree,''
{\sl arXiv preprint 1101.3656}\/
(1998).
\endref
\beginref
\parindent=20pt\item{\bref{18}}
\hldest{xyz}{}{bib18}%
Harry Kesten,
``Subdiffusive behavior of a random walk on a random cluster,''
{\sl Annales de l'Institut Henri Poincar\'e}\/
{\bf 22}
(1986),
425--487.
\endref
\beginref
\parindent=20pt\item{\bref{19}}
\hldest{xyz}{}{bib19}%
Valintin Fedorovich Kolchin,
{\sl Random Mappings}
(New York:
Optimisation Software Inc.,
1986).
\endref
\beginref
\parindent=20pt\item{\bref{20}}
\hldest{xyz}{}{bib20}%
Russell Lyons
and Yuval Peres,
{\sl Probability on Trees and Networks}
(Cambridge:
Cambridge University Press,
2017).
\endref
\beginref
\parindent=20pt\item{\bref{21}}
\hldest{xyz}{}{bib21}%
Hosam Mahmoud Mahmoud
and Mark Daniel Ward,
``Asymptotic distribution of two-protected nodes in random binary search trees,''
{\sl Applied Mathematics Letters}\/
{\bf 25}
(2012),
2218--2222.
\endref
\beginref
\parindent=20pt\item{\bref{22}}
\hldest{xyz}{}{bib22}%
Hosam Mahmoud Mahmoud
and Mark Daniel Ward,
``Asymptotic properties of protected nodes in random recursive trees,''
{\sl Journal of Applied Probability}\/
{\bf 52}
(2015),
290--297.
\endref
\beginref
\parindent=20pt\item{\bref{23}}
\hldest{xyz}{}{bib23}%
Toufik Mansour,
``Protected points in $k$-ary trees,''
{\sl Applied Mathematics Letters}\/
{\bf 24}
(2011),
478--480.
\endref
\beginref
\parindent=20pt\item{\bref{24}}
\hldest{xyz}{}{bib24}%
Boris Alexeyevich Rogozin,
``On an estimate of the concentration function,''
{\sl Theory of Probability and its Applications}\/
{\bf 6}
(1961),
94--97.
\endref
\beginref
\parindent=20pt\item{\bref{25}}
\hldest{xyz}{}{bib25}%
Benedikt Stufler,
``Local limits of large Galton--Watson trees rerooted at a random vertex,''
{\sl Annales de l'Institut Henri Poincar\'e, Probabilit\'es et Statistiques}\/
{\bf 55}
(2019),
155--183.
\endref
\beginref
\bye

%% file: fontmac.tex


\font\smallheader=cmssbx10 

\font\eightpt=cmr8
\font\ninept=cmr9


\font\mathbold=cmmib10

\font\ninerm=cmr9     \font\eightrm=cmr8   \font\sixrm=cmr6      
\font\ninei=cmmi9     \font\eighti=cmmi8   \font\sixi=cmmi6      
\font\ninesy=cmsy9    \font\eightsy=cmsy8  \font\sixsy=cmsy6     
\font\ninebf=cmbx9    \font\eightbf=cmbx8  \font\sixbf=cmbx6     
\font\ninett=cmtt9    \font\eighttt=cmtt8                        
\font\nineit=cmti9    \font\eightit=cmti8     
\font\ninesl=cmsl9    \font\eightsl=cmsl8                        

\font\tensc=cmcsc10   \font\ninesc=cmcsc9  \font\eightsc=cmcsc8  

\font\eightssq=cmssq8  \font\eightssqi=cmssqi8  


\font\tenssbx=cmssbx10 

  \font\twelvebf=cmbx12
  
\def\sc{\tensc}  \def\mc{\ninerm}

\input cyracc.def
    \font\tencyr=wncyr10   \font\ninecyr=wncyr9   \font\eightcyr=wncyr8
    \font\tencyri=wncyi10  \font\ninecyri=wncyi9  \font\eightcyri=wncyi8
    \def\cyr{\tencyr\cyracc} \def\cyri{\tencyri\cyracc}

\newskip\ttglue
\def\tenpoint{\def\rm{\fam0\tenrm}%
  \textfont0=\tenrm \scriptfont0=\sevenrm \scriptscriptfont0=\fiverm
  \textfont1=\teni  \scriptfont1=\seveni  \scriptscriptfont1=\fivei
  \textfont2=\tensy \scriptfont2=\sevensy \scriptscriptfont2=\fivesy
  \textfont3=\tenex \scriptfont3=\tenex   \scriptscriptfont3=\tenex
  \textfont\itfam=\tenit  \def\it{\fam\itfam\tenit}%
  \textfont\slfam=\tensl  \def\sl{\fam\slfam\tensl}%
  \textfont\ttfam=\tentt  \def\tt{\fam\ttfam\tentt}%
  \textfont\bffam=\tenbf  \scriptfont\bffam=\sevenbf
   \scriptscriptfont\bffam=\fivebf \def\bf{\fam\bffam\tenbf}%
  \tt \ttglue=.5em plus.25em minus.15em
  \normalbaselineskip=12pt
  \setbox\strutbox=\hbox{\vrule height8.5pt depth3.5pt width0pt}%
  \let\sc=\tensc \let\mc=\ninerm  
  \def\cyr{\tencyr\cyracc}\def\cyri{\tencyri\cyracc}
  \let\big=\tenbig  \normalbaselines\rm}

\def\ninepoint{\def\rm{\fam0\ninerm}%
\textfont0=\ninerm  \scriptfont0=\sixrm  \scriptscriptfont0=\fiverm
\textfont1=\ninei   \scriptfont1=\sixi   \scriptscriptfont1=\fivei
\textfont2=\ninesy  \scriptfont2=\sixsy  \scriptscriptfont2=\fivesy
\textfont3=\tenex   \scriptfont3=\tenex  \scriptscriptfont3=\tenex
\textfont\itfam=\nineit  \def\it{\fam\itfam\nineit}%
\textfont\slfam=\ninesl  \def\sl{\fam\slfam\ninesl}%
\textfont\ttfam=\ninett  \def\tt{\fam\ttfam\ninett}%
\textfont\bffam=\ninebf  \scriptfont\bffam=\sixbf
\scriptscriptfont\bffam=\fivebf\def\bf{\fam\bffam\ninebf}%
\tt\ttglue=.5em plus.25em minus.15em
\normalbaselineskip=11pt
\setbox\strutbox=\hbox{\vrule height8pt depth3pt width0pt}%
\let\sc=\ninesc\let\mc=\eightrm
\def\cyr{\ninecyr\cyracc}\def\cyri{\ninecyri\cyracc}
\let\big=\ninebig\normalbaselines\rm}

\def\eightpoint{\def\rm{\fam0\eightrm}%
  \textfont0=\eightrm \scriptfont0=\sixrm \scriptscriptfont0=\fiverm
  \textfont1=\eighti  \scriptfont1=\sixi  \scriptscriptfont1=\fivei
  \textfont2=\eightsy \scriptfont2=\sixsy \scriptscriptfont2=\fivesy
  \textfont3=\tenex   \scriptfont3=\tenex \scriptscriptfont3=\tenex
  \textfont\itfam=\eightit  \def\it{\fam\itfam\eightit}%
  \textfont\slfam=\eightsl  \def\sl{\fam\slfam\eightsl}%
  \textfont\ttfam=\eighttt  \def\tt{\fam\ttfam\eighttt}%
  \textfont\bffam=\eightbf  \scriptfont\bffam=\sixbf
  \normalbaselineskip=9pt
  \let\sc=\eightsc \let\mc=\sevenrm  
  \def\cyr{\eightcyr\cyracc}\def\cyri{\eightcyri\cyracc}
  \let\big=\eightbig  \normalbaselines\rm}%
\def\nospace{\nulldelimiterspace0pt\mathsurround0pt}%
\def\tenbig#1{{\hbox{$\left#1\vbox to8.5pt{}\right.\nospace$}}}%
\def\ninebig#1{{\hbox{$\textfont0=\tenrm\textfont2=\tensy
  \left#1\vbox to7.25pt{}\right.\nospace$}}}%
\def\eightbig#1{{\hbox{$\textfont0=\ninerm\textfont2=\ninesy
  \left#1\vbox to6.5pt{}\right.\nospace$}}}%

\def\nonextendedbold{
  \font\fiveb=cmb10 at 5pt
  \font\sixb=cmb10 at 6pt
  \font\sevenb=cmb10 at 7pt
  \font\eightb=cmb10 at 8pt
  \font\nineb=cmb10 at 9pt
  \font\tenb=cmb10
  \font\twelveb=cmb10 at 12pt
  \let\fivebf=\fiveb
  \let\sixbf=\sixb
  \let\sevenbf=\sevenb
  \let\eightbf=\eightb
  \let\ninebf=\nineb
  \let\tenbf=\tenb
  \let\twelvebf=\twelveb
}

\def\leftrighttop#1#2{
  \headline{\ifnum\pageno=1\hfil\else{\ninept #1 \hfil #2}\fi}
}

\def\firstnopagenum{
  \footline{\ifnum\pageno=1 \hfil \else \hfil{\rm \number\pageno}\hfil\fi}
}

\def\maketitle#1#2#3#4{
  \centerline {\titlefont #1}
  \medskip
  \centerline {\eightpt #2}
  \medskip
  \centerline {\tensc #3}
  \medskip
  \centerline {\tensc #4}
  \bigskip
}


\outer\def\floattext#1 #2. #3\par{
  $$
  \vbox{
    \hsize #1 true in
    \noindent{\bf #2.}\enskip #3
  }
  $$
}


\def\lsection#1\par{
  \bigskip\vskip\parskip
  \leftline{\sectionfont#1}\nobreak\medskip\noindent
}

\def\csection#1\par{
  \bigskip\vskip\parskip
  \centerline{\sectionfont#1}\nobreak\medskip\noindent
}

\def\rsection#1\par{
  \bigskip\vskip\parskip
  \rightline{\sectionfont#1}\nobreak\medskip\noindent
}
\def\section{\lsection}

\def\boldlabel#1. {\noindent{\bf #1.\enspace}}
\def\subsection#1. {\medskip\noindent{\bf #1.\enspace}}

\def\caption Fig. #1. #2.{\ninepoint{\bf Fig.\ #1.}\enspace#2.}



\font\tenfrak=eufm10
\font\sevenfrak=eufm7
\font\fivefrak=eufm5
\newfam\frakfam
\textfont\frakfam=\tenfrak
\scriptfont\frakfam=\sevenfrak
\scriptscriptfont\frakfam=\fivefrak

\def\janksc#1#2 {#1{\eightpt#2}}
\def\jankscsp#1#2 {#1{\eightpt#2}\ }
\def\scproclaim#1.#2\par{\noindent\jankscsp #1.\enspace{\it#2\par}}

\def\ref#1{[#1]}

\def\quote{
  \begingroup
    \baselineskip 10pt
    \parfillskip 0pt
    \interlinepenalty 10000
    \leftskip 0pt plus 40pc minus \parindent
    \let\rm=\quoterm\let\sl=\quotesl\everypar{\sl}
    \obeylines
}
\def\author#1(#2){\nobreak\smallskip\rm--- \rm#1\unskip\enspace(#2)\par\endgroup}

\def\titlefont{\twelvebf}
\def\sectionfont{\tenssbx}
\def\quoterm{\eightssq}
\def\quotesl{\eightssqi}


\tenpoint

%% file: mathmac.tex


\def\xskip{\hskip 7pt plus 3pt minus 4pt}

\def\proof{\medbreak\noindent{\it Proof.}\xskip\ignorespaces}

\def\slug{\quad\hbox{\kern1.5pt\vrule width2.5pt height6pt depth1.5pt\kern1.5pt}\medskip}
\def\noskipslug{\quad\hbox{\kern1.5pt\vrule width2.5pt height6pt depth1.5pt\kern1.5pt}}

\newdimen\algindent
\newif\ifitempar \itempartrue 
\def\algindentset#1{\setbox0\hbox{{\bf #1.\kern.25em}}\algindent=\wd0\relax}
\def\algbegin #1 #2{\algindentset{#21}\alg #1 #2} 
\def\aalgbegin #1 #2{\algindentset{#211}\alg #1 #2} 
\def\alg#1(#2). {\medbreak 
  \noindent{\bf#1}({\it#2\/}).\xskip\ignorespaces}
\def\algstep#1.{\ifitempar\smallskip\noindent\else\itempartrue
  \hskip-\parindent\fi
  \hbox to\algindent{\bf\hfil #1.\kern.25em}%
  \hangindent=\algindent\hangafter=1\ignorespaces}




\def\pr{\hbox{\bf P}}
\def\ex{\hbox{\bf E}}
\def\var{\hbox{\bf V}}



\def\eps{\epsilon}

\newcount\thmcount  
\thmcount=1
\newcount\sectcount  
\sectcount=1
\newcount\figcount  
\figcount=1
\newcount\eqcount  
\eqcount=1

\def\oldno#1{\eqno({\oldstyle#1})}
\def\refeq#1{({\oldstyle#1})}
\def\adveq{\oldno{\the\eqcount}\global\advance\eqcount by 1}  
\def\advthm{\the\thmcount\global\advance \thmcount by 1}

\def\caption#1{\centerline{\ninepoint{\bf Fig.~\the\figcount\global\advance\figcount by 1.\enspace}#1}}

\def\advsect{\section\the\sectcount\global\advance\sectcount by 1. }

\outer\def\parenproclaim #1 (#2).#3\par{\medbreak
  \noindent{\bf #1}\enspace\rm({\it #2\/}).\nobreak\ignorespaces{\sl #3\par}
  \ifdim\lastskip<\medskipamount \removelastskip\penalty55\medskip\fi}
